\documentclass[12pt,a4paper]{article}

\usepackage[english]{babel}
\usepackage{epsfig}
\usepackage{amsmath}
\usepackage{amssymb}
\usepackage{amsbsy}
\usepackage{amsfonts}
\usepackage{mathrsfs} 
\usepackage{array}
\usepackage{verbatim}
\usepackage{graphicx}
\usepackage{color}
\usepackage{tabularx}
\usepackage{bm}
\usepackage{booktabs}
\usepackage[font=small]{caption}
\usepackage{subfig}
\usepackage{stmaryrd}
\usepackage{textcomp}

\usepackage{booktabs}
\usepackage{array}
\usepackage{float}

\usepackage{amsmath}
\usepackage{amssymb}
\usepackage{amsthm}

\usepackage{graphicx}
\usepackage{IEEEtrantools}
\usepackage{mdframed} 

\usepackage{color}
\definecolor{light}{gray}{0.50}
\definecolor{heavy}{gray}{0.35}
\definecolor{black}{gray}{0.0}
\definecolor{dgreen}{rgb}{0.0,0.7,0}
\definecolor{dred}{rgb}{0.9959,0,0}
\definecolor{green}{rgb}{0.0,0.99599,0.0}
\definecolor{purple}{rgb}{0.6,0.0,0.4}

\theoremstyle{definition}
\newtheorem{algorithm}{Algorithm}

\let\dss=\displaystyle

\renewcommand{\vector}[1]{\relax\ifmmode\mathchoice
{\mbox{\boldmath$\displaystyle#1$}}
{\mbox{\boldmath$\displaystyle#1$}}
{\mbox{\boldmath$\scriptstyle#1$}}
{\mbox{\boldmath$\scriptscriptstyle#1$}}\else
\hbox{\boldmath$\textstyle#1$}\fi}

\newcommand{\eq}[1]{(\ref{#1})}

\newcommand{\vect}[1]{{\mathbf{#1}}}

\newcommand{\Pbar}{\overline{P}}

\newcommand{\vk}{\vect{k}}
\newcommand{\vu}{\vect{u}}
\newcommand{\vv}{\vect{v}}
\newcommand{\vx}{\vect{x}}

\newcommand{\vV}{\vect{V}}

\newcommand{\halff}{\frac{1}{2}}
\newcommand{\half}{1/2}

\newcommand{\dt}{\Delta t}
\newcommand{\dz}{\Delta z}
\newcommand{\chibar}{\overline{\chi}}

\newcommand{\chiprime}{{\chi'}}

\newcommand{\Thetabar}{\overline{\Theta}}
\newcommand{\Thetatilde}{{\widetilde \Theta}}
\newcommand{\pibar}{\overline{\pi}}

\newcommand{\piprime}{\pi'}
\newcommand{\piprimehat}{\widehat{\pi}'}

\newcommand{\Uhat}{\widehat{U}}
\newcommand{\Uhathat}{\widehat{\widehat{U}}}

\newcommand{\bigoh}[1]{\mathcal{O}\left(#1\right)}

\newcommand{\Sol}{\mathcal{U}}

\newcommand{\rfr}[1]{#1_{\text{ref}}}

\newcommand{\nablatilde}{{\widetilde\nabla}}

\newcommand{\dx}{{\Delta x}}
\newcommand{\dy}{{\Delta y}}

\newcommand{\ahydro}{\alpha_{w}}
\newcommand{\apsinc}{\alpha_{P}}

 \graphicspath{{figures/}}

\begin{document}

\title{A semi-implicit compressible model for atmospheric flows with seamless access to soundproof and hydrostatic dynamics}
\author{Tommaso Benacchio $^{(1)}$, Rupert Klein $^{(2)}$}
\maketitle

\begin{center}
 
{\small
$^{(2)}$ MOX -- Modelling and Scientific Computing, \\
Dipartimento di Matematica, Politecnico di Milano \\
Piazza Leonardo da Vinci 32, 20133 Milano, Italy\\
{\tt tommaso.benacchio@polimi.it}
}
\vskip 0.5cm
{\small
$^{(2)}$ FB Mathematik \& Informatik, \\
Freie Universit\"at Berlin, \\
Arnimallee 6, 14195 Berlin, Germany\\
{\tt rupert.klein@math.fu-berlin.de}
}

\end{center}

\date{}

\noindent
{\bf Keywords}:   semi-implicit models; finite volume methods; hyperbolic equations; compressible flow;
soundproof models; hydrostatic models.

\vspace*{0.5cm}

\noindent
{\bf AMS Subject Classification}:   65M08, 65Z99, 76B55, 76M12,  76U05, 86A10

\vspace*{0.5cm}

\pagebreak

\abstract{We introduce a second-order numerical scheme for compressible atmospheric motions at small to planetary scales. The collocated finite volume method treats the advection of mass, momentum, and mass-weighted potential temperature in conservation form while relying on Exner pressure for the pressure gradient term. It discretises the rotating compressible equations by evolving full variables rather than perturbations around a background state, and operates with time steps constrained by the advection speed only. Perturbation variables are only used as auxiliary quantities in the formulation of the elliptic problem. Borrowing ideas on forward-in-time differencing, the algorithm reframes the authors' previously proposed schemes into a sequence of implicit midpoint, advection, and implicit trapezoidal steps that allows for a time integration unconstrained by the internal gravity wave speed. Compared with existing approaches, results on a range of benchmarks of nonhydrostatic- and hydrostatic-scale dynamics are competitive. The test suite includes a new planetary-scale inertia-gravity wave test highlighting the properties of the scheme and its large time step capabilities. In the hydrostatic-scale cases the model is run in pseudo-incompressible and hydrostatic mode with simple switching within a uniform discretization framework. The differences with the compressible runs return expected relative magnitudes. By providing seamless access to soundproof and hydrostatic dynamics, the developments represent a necessary step towards an all-scale blended multimodel solver.}

\pagebreak

\section{Introduction}
\label{sec:Intro}


\subsection{Motivation: Blending of full and reduced dynamical flow models}

Atmospheric dynamics features a variety of scale-dependent motions which have been analytically described by scale analysis and asymptotics \cite{Klein2010,Pedlosky1992}. Reduced dynamical models emerging from the full compressible flow equations through generally singular asymptotic limits capture the essence of the phenomena of interest and reveal which effects are important -- and which ones less so -- for their description. Relevant examples include the anelastic and pseudo-incompressible models, the quasi-geostrophic and semi-geostrophic models, and the hydrostatic primitive model equations \cite{Bannon1996,CullenMaroofi2003,Durran1989,HoskinsBretherton1972,Klein2010,LippsHemler1982,Pedlosky1992}. 

In \cite{Cullen2007}, the author argues that compressible atmospheric flow solvers should accurately reproduce the effective dynamics encoded by such reduced dynamical models with no degradation of solution quality as the respective limit regime is approached. Related numerical methods are known as \textit{asymptotic preserving} or \textit{asymptotically adaptive} schemes in the numerics literature, see \cite{Jin2012, KleinEtal2001} for references. If a scheme is designed such that it not only solves the compressible equations close to some limit regimes with the required accuracy but that it can also solve the limiting model equations when the respective asymptotic parameter is set to zero, this opens avenues to interesting applications and investigations.

Implementations of different model equations normally use different numerical methods to represent identical terms. For example, in a comparison of a compressible model with a pseudo-incompressible model, the former might discretize advection with a semi-Lagrangian scheme, while the latter uses a higher-order upwind finite volume formulation. In this case, differences in model results cannot be uniquely attributed to the differences in the underlying equations but may as well be influenced by the use of different advection schemes (see \cite{BenacchioEtAl2014,SmolarkiewiczDoernbrack2007} for further examples).

Using a numerical method for the compressible equations that defaults to soundproof dynamics for vanishing Mach number, in \cite{BenacchioEtAl2014} the authors suggested an application in the context of well-balanced data assimilation. They implement a \textit{blended} scheme that can be tuned to solve any one of a continuous family of equations that interpolate between the compressible and pseudo-incompressible models, and use this feature to filter unwanted acoustic noise from some given or assimilated initial data. To properly capture a compressible flow situation with unknown balanced initial pressure distribution, they operate the scheme for some initial time steps in its pseudo-incompressible mode and then relax the model blending parameter towards its compressible mode over a few more steps. In this fashion, the pseudo-incompressible steps serve to find a balanced pressure field compatible with the velocity and potential temperature initial data, and the subsequent compressible flow simulation is essentially acoustics-free.

Continuing this line of development, we describe in this paper a semi-implicit scheme that allows us to access the compressible, pseudo-incompressible, and hydrostatic models within one and the same finite volume framework. 


\subsection{Related numerical schemes in the literature}

A significant challenge in the dynamical description and forecast of weather and climate lies in the inherently multiscale nature of atmospheric flows. Driven by stratification and rotation, physical processes arise around a large-scale state of horizontal geostrophic, vertical hydrostatic balance. The compressible Euler equations are deemed the most comprehensive model to describe the principal fluid dynamical features of the system before parameterizations of unresolved processes are added. On the one hand, these equations allow for buoyancy-driven internal gravity wave and pressure-driven sound wave adjustments.  On the other hand, meteorologically relevant features such as cyclones and anticyclones in the midlatitudes involve motions much slower than the sound speed, thus forcing numerical stiffness into discretizations of the compressible model in the low Mach number regime. As a result, most if not all numerical schemes used in operational weather forecasting employ varying degrees of implicitness or multiple time stepping that enable stable runs with long time step sizes unconstrained by sound speed (see, e.g. the reviews \cite{MarrasEtAl2016, MengaldoEtAl2018} and references therein for a list). Typically, semi-implicit approaches integrate advective transport explicitly, then build an elliptic problem for the pressure variable by combining the equations of the discrete system. The solution of the problem yields updates that are then replaced into the other variables. 

Examples of operational dynamical cores using semi-implicit time-integrations strategies are the ECMWF\footnote{European Centre for Medium-Range Weather Forecasts}'s IFS \cite{Hortal2002}, that discretizes the hydrostatic primitive equations, and the UK Met Office's ENDGame \cite{BenacchioWood2016, WoodEtAl2013}. In particular, ENDGame  uses a double-loop structure in the implicit solver entailing four solves per time step in its operational incarnation, a strategy carried over in recent developments \cite{MelvinEtAl2018}, and allowing non-operational configurations to run stably and second-order accurately without additional numerical damping (for operational forecasts, a small amount of off-centering is usually employed for safety reasons). By contrast, many other semi-implicit or time-split explicit discretizations resort to off-centering, divergence damping \cite{BryanFritsch2002}, or otherwise artificial diffusion in order to quell numerical instabilities. In non-operational research, the authors of \cite{DumbserEtAl2018}, among others, present buoyancy- and acoustic-implicit second-order finite volume discretisations on staggered grids.

In order to simplify the formulation of the semi-implicit method, the equation set is often cast in terms of perturbations around a hydrostatically balanced reference state, see, e.g., \cite{RestelliGiraldo2009,SmolarkiewiczEtAl2019,SmolarkiewiczEtAl2014}. However, as noted in \cite{WellerShahrokhi2014}, whose model does not use pertubations, large deviations from the background state may question the assumptions underpinning the resulting system. In \cite{WoodEtAl2013} and \cite{MelvinEtAl2018}, the authors use the model state computed at the previous time step as evolving background profile.

To address the efficiency issues caused by spectral transforms in IFS at increasing global resolutions, a finite volume discretization is also used in FVM, the potential next-generation ECMWF dynamical core \cite{KuehnleinEtAl2019}. The time integration algorithm in FVM is built upon extensive earlier experience with the EULAG model and the MPDATA advection scheme. Through appropriate correction of a first-order upwind discretization, a system is constructed that encompasses transport and implicit dynamics in an elegant theoretical framework \cite{SmolarkiewiczEtAl2016, SmolarkiewiczEtAl2014} and references therein. The approach, which in its default configuration relies on time extrapolation of advecting velocities and subtraction of reference states, also contains soundproof analytical systems as subcases and has shown excellent performance in integrating atmospheric flows at all scales without instabilities. However, their transition from compressible to soundproof discretizations is not seamless in the sense of the present work, since the structure of their implicit pressure equations substantially differs from one model to the other. Similarly to the present approach, an optional variant of their scheme avoids extrapolations in time from earlier time levels.

Drawing on the finite volume framework for soundproof model equations in \cite{KleinTCFD2009}, the authors of \cite{Benacchio2014, BenacchioEtAl2014} devised a numerical scheme for the compressible Euler equations to simulate small- to mesoscale atmospheric motions, using a time step unconstrained by the speed of acoustic waves within the abovementioned soundproof-compatible switchable multimodel formulation. The underlying theoretical framework was extended in \cite{KleinBenacchio2016} to incorporate the hydrostatic primitive equations and the anelastic, quasi-hydrostatic system \cite{ArakawaKonor2009} with the introduction of a second blending parameter.

A major hurdle towards joining the numerical scheme of \cite{BenacchioEtAl2014} with the theoretical setup of \cite{KleinBenacchio2016} is the former's time step dependency on the speed of internal gravity waves, a severe constraint on the applicability of the numerical method to large-scale tests. The present study addresses this fundamental shortcoming. 

\subsection{Contribution}

By reframing the schemes in \cite{KleinTCFD2009} and \cite{BenacchioEtAl2014} as a two-stage-implicit plus transport system, this paper proposes a discretisation that:
\begin{itemize}
\item Evolves the compressible equations with rotation in terms of full variables, using auxiliary perturbation variables only in formulating the buoyancy-implicit elliptic problem;
\item Has built-in conservation of mass and mass-weighted potential temperature, and is second-order accurate in all components, without artificial damping mechanisms;
\item Uses a time step constrained only by the underlying advection speed;
\item Works with a node-based implicit pressure equation only, thereby avoiding the usual cell-centered MAC-projection (see \cite{AlmgrenEtAl1998} and references therein); 
\item Can be operated in the soundproof and hydrostatic modes without modifying the numerics;
\item Constitutes a basis for a multiscale formulation with access to hydrostasy and geostrophy.
\end{itemize}

The method uses an explicit second-order MUSCL scheme for advection, while the pressure and momentum equations are stably integrated by solving two elliptic problems embedded in the implicit midpoint and implicit trapezoidal stages. The scheme is validated against two-dimensional Cartesian benchmarks of nonhydrostatic and hydrostatic dynamics. Simulations of inertia-gravity wave tests at large scale and with rotation show competitive performance with existing approaches already at relatively coarse resolutions. In particular, a new planetary-scale extension of the hydrostatic-scale test of \cite{SkamarockKlemp1994} showcases the large time step capabilities of the present scheme. For the large-scale tests, we run the model in pseudo-incompressible mode and hydrostatic mode and analyse the difference with the compressible simulation. As expected from theoretical normal mode analyses \cite{DaviesEtAl2003, Dukowicz2013} (though see also \cite{KleinEtAl2010} for a discussion on regime of validity of soundproof models), the compressible/hydrostatic discrepancy shrinks with smaller vertical-to horizontal domain size aspect ratios, while the compressible/pseudo-imcompressible discrepancy grows with larger scales.

The paper is organized as follows. Section \ref{sec:GoverningEquations} contains the governing equations that are discretized with the methodology summarised in section \ref{sec:TimeDiscretizationSummary} and detailed in section~\ref{sec:DiscretizationDetails}. Section \ref{sec:Results} documents the performance of the code on the abovementioned tests. Results are discussed and conclusions drawn in section \ref{sec:Conclusions}. 


\section{Governing equations}
\label{sec:GoverningEquations}

The governing equations for adiabatic compressible flow of an inert ideal gas 
with constant specific heat capacities under the influence of gravity and in 
a rotating coordinate system corresponding to a tangent plane approximation
may be written as
\begin{IEEEeqnarray}{rCl}\label{eq:CompressibleEuler}
\dss \rho_t + \nabla_\parallel\cdot(\rho \vu) + (\rho w)_z
  & = 
    & \dss 0
      \IEEEyesnumber\IEEEyessubnumber*\label{eq:EulerMass}\\[5pt]
\dss (\rho\vu)_t + \nabla_\parallel\cdot(\rho \vu\circ\vu) + (\rho w \vu)_z 
  & = 
    & \dss - \left[ c_p  P \nabla_\parallel \pi + f(y) \vk \times \rho\vu \right]
      \label{eq:EulerHorMom}\\[5pt]
\dss (\rho w)_t + \nabla_\parallel\cdot(\rho \vu w) + (\rho w^2)_z 
  & = 
    & \dss - \left(  c_p P \pi_z + \rho g \right)
      \label{eq:EulerVerMom}\\[5pt]
\dss P_t + \nabla_\parallel\cdot(P\vu) + (Pw)_z
  & = 
    & \dss 0\,.
    \label{eq:EulerPressure}
\end{IEEEeqnarray}
Here $\rho$ is the density, $\vu = (u,v)$ and $w$ are the horizontal and vertical 
components of the flow velocity,  
\begin{equation}\label{eq:EOSpiP}
\pi = \left(\frac{p}{\rfr{p}}\right)^{\frac{R}{c_p}}
\qquad\text{and}\qquad
P = \frac{\rfr{p}}{R} \left(\frac{p}{\rfr{p}}\right)^{\frac{c_v}{c_p}} \equiv \rho\Theta
\end{equation}
are the Exner pressure and the mass-weighted potential temperature, with $\rfr{p}$ a suitable reference pressure, $R$ the gas constant and $c_p$ and 
$c_v = c_p - R$ the 
specific heat capacities at constant pressure and constant volume. Furthermore, $g$ is the acceleration of gravity (taken as constant), $f(y) = f_0 + \beta y$ the local Coriolis parameter in 
the $\beta$-plane with constant $f_0$ and $\beta$, $\vk$ the vertical 
unit vector, and $\times$ the cross product. Subscripts as in 
$U_x \equiv \partial_x U := \partial U/ \partial x$ denote partial derivatives with respect 
to the first coordinate of a Cartesian $(x,y,z)$ coordinate system or time $t$, and 
$\nabla_\parallel = (\partial_x, \partial_y, 0)$ subsumes the horizontal derivatives.

Given \eq{eq:EulerMass} and \eq{eq:EulerPressure}, the potential temperature
$\Theta = P/\rho$ satisfies the usual advection equation
\begin{equation}
\Theta_t + \vu\cdot\nabla_\parallel \Theta + w \Theta_z = 0\,.
\end{equation}


\section{Compact description of the time integration scheme}
\label{sec:TimeDiscretizationSummary}

In this section we describe the main structural features of the discretization, which evolves and joins aspects of the models in \cite{KleinTCFD2009,BenacchioEtAl2014}, and borrows key ideas from the \emph{forward-in-time} integration strategy \cite{SmolarkiewiczMargolin1997, SmolarkiewiczMargolin1993} in realizing the implicit discretization of the gravity term. 


\subsection{Reformulation of the governing equations}
\label{ssec:Reformulation}


\subsubsection{Evolution of the primary variables}
\label{sssec:PrimaryVariables}

The primary unknowns advanced in time by the present scheme are the same as
in \eq{eq:CompressibleEuler}, i.e., $(\rho, \rho\vu, \rho w, P)$. 
Introducing a seamless blended discretization of the compressible Euler and
pseudo-incompressible equations \cite{Durran1989} and following
\cite{KleinTCFD2009,KleinEtAl2010}, in \cite{BenacchioEtAl2014} the authors observed that 
the pseudo-incompressible model is obtained from the compressible equations in 
\eq{eq:CompressibleEuler} by simply dropping the time derivative of 
$P = \rho\Theta$ from \eq{eq:EulerPressure}. To take advantage of this simple 
structural model relationship in constructing a blended scheme that can be
tuned seamlessly from solving the full compressible model equations
to solving the pseudo-incompressible model equations, they introduced the inverse of the 
potential temperature,
\begin{equation}
\chi = 1/\Theta\,,
\end{equation}
and interpreted the mass balance \eq{eq:EulerMass} as a transport equation for $\chi$, 
\begin{equation}\label{eq:chiI}
\rho_t + \nabla_\parallel\cdot(\rho \vu) + (\rho w)_z = 
(P\chi)_t + \nabla_\parallel\cdot(P\chi \vu) + (P\chi w)_z = 0\,,
\end{equation}
in which the field $(P\vv)$ takes the role of an advecting flux. 
Using this interpretation consistently throughout the equation system, and introducing two blending parameters, $\ahydro$ and $\apsinc$, for the non-hydrostatic/hydrostatic and compressible/pseudo-incompressible transitions, one obtains
\begin{IEEEeqnarray}{rCrCl}\label{eq:EulerP}
\dss \rho_t 
  & + 
    & \dss \nabla_\parallel\cdot(P\vu\, \chi) + (Pw\, \chi)_z \hfil
      & = 
        & \dss 0
          \IEEEyesnumber\IEEEyessubnumber*\label{eq:EulerPMass}\\[5pt]
\dss (\rho\vu)_t 
  & + 
    & \dss \nabla_\parallel\cdot(P\vu\circ\chi\vu) + (Pw\, \chi\vu)_z  \hfil
      & = 
        & \dss - \left[ c_p P\nabla_\parallel \pi \right.
          \label{eq:EulerPHorMom}\\[5pt]
       & & & &  \left. + f(y) \vk\times\rho\vu\right]\nonumber\\
\dss \ahydro \Bigl[(\rho w)_t 
  & + 
    & \dss \nabla_\parallel\cdot(P \vu\, \chi w) + (Pw\, \chi w)_z\Bigr] \hfil
      & = 
        & \dss - \left( c_p P \pi_z + \rho g\right)
          \label{eq:EulerPVerMom}\\[5pt]
\apsinc\, P_t
  &  +
    & \dss \dss \nabla_\parallel\cdot(P\vu)  + (Pw)_z  \hfil
      & = 
        & \dss 0\,.
        \label{eq:EulerPP}
\end{IEEEeqnarray}
System \eq{eq:EulerP} is the analytical formulation used in this paper, and facilitates the extension of the blending of \cite{BenacchioEtAl2014} to hydrostasy along the lines of the theory described in \cite{KleinBenacchio2016}. The quasi-geostrophic case will be addressed in forthcoming work.


\subsubsection{Auxiliary perturbation variables and their evolution equations}
\label{sssec:AuxPerturbationVariables}

A crucial ingredient of any numerical scheme implicit with respect to the 
effects of compressibility, buoyancy, and Earth rotation, is that it has separate 
access to the large-scale mean background stratifications of pressure and 
potential temperature, or its inverse, and to their local perturbations. 
Thus we split the Exner pressure $\pi$ and inverse potential temperature $\chi$ into
\begin{equation}\label{eq:PerturbationVariables}
\pi(t,\vx,z) = \piprime(t,\vx,z) + \pibar(z)
\qquad\text{and}\qquad
\chi(t,\vx,z) = \chiprime(t,\vx,z) + \chibar(z)\, ,
\end{equation}
with the hydrostatically balanced background variables satisfying
\begin{equation}\label{eq:BackgroundHydrostatics}
\frac{d\pibar}{dz} = - \frac{g}{c_p} \chibar
\qquad\text{and}\qquad
\pibar(0) = 1\, .
\end{equation}
Since, for the compressible case, $P$ can be expressed as a function of $\pi$ alone according to
\eq{eq:EOSpiP}, and since $\pibar$ is time independent across a time step, 
the perturbation Exner pressure satisfies
\begin{equation}\label{eq:EulerPiPrime}
\apsinc \, \left(\frac{\partial P}{\partial \pi}\right) \piprime_t
= 
- \nabla\cdot \left[P(\pi)\vv\right]\,,
\end{equation}
which is a direct consequence of \eq{eq:EulerPP}.
In turn, the perturbation form of the mass balance serves as the evolution equation
for $\chiprime$, i.e.,
\begin{IEEEeqnarray}{rCrCl}\label{eq:EulerPChiPrime}
\dss (P\chiprime)_t 
  & + 
    & \dss \nabla_\parallel\cdot(P\vu\, \chiprime) + (Pw\, \chiprime)_z \hfil
      & = 
        & \dss -\left[\nabla_\parallel\cdot(P\vu\, \chibar) + (Pw\, \chibar)_z\right]\,.
        \label{eq:ChiPrimeEqn}
\end{IEEEeqnarray}

Auxiliary discretizations of \eq{eq:EulerPiPrime} and 
\eq{eq:EulerPChiPrime} will be used in constructing a numerical scheme
for the full variable form of the governing equations in \eq{eq:EulerP}
that is stable for time steps limited only by the advection Courant 
number. After completion of a time step, the perturbation variables
are synchronized with the full variables based on the definitions
in \eq{eq:PerturbationVariables} and \eq{eq:BackgroundHydrostatics}. We remark that this is a fundamental feature of the present scheme, shared with the staggered grid scheme by \cite{WellerShahrokhi2014}. To the best of our knowledge, other models for atmospheric flows use the perturbation variables as prognostic quantities throughout.

In the sequel, borrowing notation from \cite{SmolarkiewiczEtAl2014},
we introduce
\begin{equation}\label{eq:PsiDefinition}
\Psi = (\chi, \chi\vu, \chi w, \chiprime)
\end{equation}
and subsume the primary equations in \eq{eq:EulerP} and the auxiliary 
equation for $\chiprime$ in \eq{eq:ChiPrimeEqn} as 
\begin{IEEEeqnarray}{rCl}\label{eq:PBasedAdvection}
(P\Psi)_t + \mathcal{A}(\Psi; P\vv) 
  & = 
    & Q(\Psi; P)
      \IEEEyesnumber\IEEEyessubnumber*\label{eq:EulerPCompactPsi}\\
\apsinc \,P_t + \nabla\cdot(P\vv)
  & =
    & 0\,.
    \label{eq:EulerPCompactP}
\end{IEEEeqnarray}
Note that the $\pi'$ equation in \eq{eq:EulerPiPrime} is 
equivalent to \eq{eq:EulerPCompactP} and thus it is not listed separately,
although it will be used in an auxiliary step in the design of a stable
discretization of \eq{eq:EulerPCompactP}.


\subsection{Semi-implicit time discretization}
\label{ssec:TimeDiscretizationOverview}


\subsubsection{Implicit midpoint pressure update and advective fluxes}
\label{sssec:AdvectiveFluxes}

In the first step of the scheme, we determine advective fluxes 
at the half-time level, $(P\vv)^{n+\half}$, which for $\apsinc = 1$ 
immediately yield the update of the internal energy variable, $P$, through
\begin{equation}\label{eq:PUpdate}
\apsinc\left(P^{n+1} - P^{n}\right)
= - \dt \,\nablatilde\cdot(P\vv)^{n+\half}\,,
\end{equation}
where $\nablatilde \cdot$ is the discrete approximation of the divergence.
In contrast, for $\apsinc = 0$ this equation represents the pseudo-incrompressible
divergence constraint.

Note that in the compressible case this update corresponds to a time 
discretization of the $P$-equation
using the \emph{implicit midpoint rule}. We recall here for future reference that
an implementation of the implicit midpoint rule can be achieved by first applying
a half time step based on the implicit Euler scheme followed by another half time
step based on the explicit Euler method \cite{HairerEtAl2006}.

To maintain second-order accuracy of the overall scheme, a first-order accurate 
time integration from the last completed time step at $t^n$ is sufficient for 
generating the half time level fluxes $(P\vv)^{n+\half}$. This becomes transparent through a truncation error analysis for any equation of 
the form $\dot y = R(y,t)$. First we observe that
\begin{equation}\label{eq:FirstOrderHTLUpdateI}
\frac{y(t^{n+1})-y(t^n)}{\dt} = \dot y\left(t^{n+\half}\right) + \bigoh{\dt^2}
\end{equation}
by straightforward Taylor expansion. Then, for any first-order approximation, 
say $R^{n+\half}$, to the right hand side at the half time level we have  
\begin{equation}\label{eq:FirstOrderHTLUpdateII}
\dot y\left(t^{n+\half}\right) 
= R\left[y\left(t^{n+\half}\right)\right] 
= R\left[y(t^n) + \frac{\Delta t}{2} \dot y(t^n) + \bigoh{\Delta t^2}\rule{0pt}{12pt}\right]
= R^{n+\half} + \bigoh{\dt^2}\,,
\end{equation}
where $R^{n+\half} = R\left[y(t^n) + (\Delta t/2) \dot y(t^n)\rule{0pt}{10pt}\right]$ is the right hand side evaluated at a state that is lifted from $t^n$ to $t^{n+\half}$ just by a first-order method. Re-inserting into \eq{eq:FirstOrderHTLUpdateI} we find indeed
\begin{equation}
\frac{y(t^{n+1})-y(t^n)}{\dt} = R^{n+\half} + \bigoh{\dt^2}\, .
\end{equation}

In order to achieve stability for large time steps, only limited by the advection Courant number,
we invoke standard splitting into advective and non-advective terms in 
\eq{eq:EulerP}, \eq{eq:EulerPChiPrime} for the prediction of $(P\vv)^{n+\half}$, 
with explicit advection and 
linearly implicit treatment of the right hand sides. 
Thus we first advance the scalars from \eq{eq:PsiDefinition} by half an advection 
time step using advective fluxes computed at the old time level, 
\begin{IEEEeqnarray}{rCl}\label{eq:HalfTimePredictorAdvection}
\dss (P\Psi)^{\#} 
  & = 
    & \dss \mathcal{A}_{1\text{st}}^{\frac{\dt}{2}}\left(\Psi^{n}; (P\vv)^{n}\right)
      \IEEEyesnumber\IEEEyessubnumber*\label{eq:HalfTimePredictorAdvectionA}\\
\dss P^{\#} 
  & = 
    & \dss P^{n} - \frac{\dt}{2} \nablatilde\cdot(P\vv)^{n}\,.
    \label{eq:HalfTimePredictorAdvectionII}
\end{IEEEeqnarray}
Here $\mathcal{A}_{1\text{st}}^{\dt}$ denotes an at least first-order accurate version of our advection scheme for the $\Psi$-variables given the advecting fluxes $(P\vv)^{n}$, see section~\ref{ssec:Advection} for details.  In the pseudo-incompressible case the discretization guarantees that $(P\vv)^n$ is discretely divergence free as shown below, so that $P^{\#} = P^n$ and the $\apsinc$ parameter need not be explicitly noted in \eq{eq:HalfTimePredictorAdvectionII}.

Next, the half time level fluxes $(P\vv)^{n+\half}$ are obtained via the implicit Euler discretization
of a second split system that only involves the right hand sides of \eq{eq:EulerP} (see section~\ref{ssec:SemiImplicitForcing} below for details), 
\begin{IEEEeqnarray}{rCl}\label{eq:HalfTimePredictorFluxCorrection}
\dss (P\Psi)^{n+\half} 
  & = 
    & \dss (P\Psi)^{\#} + \frac{\dt}{2} Q\left(\Psi^{n+1/2}; P^{n+\half}\right)\,,
      \IEEEyesnumber\IEEEyessubnumber*\label{eq:HalfTimePredictorRHS}\\
\dss \apsinc\, P^{n+\half} 
  & = 
    & \dss \apsinc\, P^{n} - \frac{\dt}{2} \nabla\cdot(P\vv)^{n+\half}\,.
      \label{eq:HalfTimePredictorDivControl}
\end{IEEEeqnarray}
We note that for $\apsinc = 1$ \eq{eq:HalfTimePredictorDivControl} corresponds to the
implicit Euler update of $P$ to the half time level, \emph{i.e.}, to the first step of 
our implementation of the implicit midpoint rule for this variable.
Furthermore, as in \cite{BenacchioEtAl2014}, in this step the relation between $P$, which is being updated by the flux divergence, and $\pi$, whose gradient is part of the momentum forcing terms, is approximated through a 
linearization of the equations of state \eq{eq:EOSpiP},  
\begin{IEEEeqnarray}{rCl}\label{eq:HalfTimePredictorPLinearization}
\dss P^{n+\half} 
  & = 
    & \dss P^{n} 
      + \left(\frac{\partial P}{\partial \pi}\right)^{\#} 
        \left(\pi^{n+\half} - \pi^{n}\right)\,.
\end{IEEEeqnarray}
With this linearization, this implicit Euler step involves a single linear elliptic 
solve for $\pi^{n+\half}$. Optionally, an outer iteration 
of the linearly implicit step can be invoked to guarantee consistency with the 
equation of state for $P(\pi)$ up to a given tolerance. 

These preliminary calculations serve to provide the fluxes $(P\vv)^{n+\half}$ later
needed both for the final explicit Euler update of $P$ to the full
time level $t^{n+1}$ and for the advection of the vector of specific variables $\Psi$ 
from \eq{eq:PsiDefinition} as part of the overall time stepping algorithm, 
see \eq{eq:AdvectionStep} below.

For $\apsinc = 0$ the $P$ equation reduces to the pseudo-incompressible
divergence constraint, and $P$ and the Exner pressure $\pi$ decouple. While $P \equiv \Pbar(z)$
remains constant in time in this case, increments of $\pi$ correspond to the elliptic pressure field
that guarantees compliance of the velocity with the divergence constraint.
 

\subsubsection{Implicit trapezoidal rule along explicit Lagrangian 
paths for advected quantities}
\label{sssec:FullTimeStep}

Given the advective fluxes, $(P\vv)^{n+\half}$, the full second-order semi-implicit time step for the evolution equation of the advected scalars, $\Psi$, reads
\begin{IEEEeqnarray}{rCl}\label{eq:TimeIntegrator}
\dss (P\Psi)^{*} 
  & = 
    & \dss (P\Psi)^{n} + \frac{\dt}{2} Q\left(\Psi^n; P^n\right)
      \IEEEyesnumber\IEEEyessubnumber*\label{eq:ExplicitEulerStep}\\
\dss (P\Psi)^{**} 
  & = 
    & \dss \mathcal{A}_{2\text{\tiny nd}}^{\dt}\left(\Psi^*; (P\vv)^{n+\half}\right)
      \label{eq:AdvectionStep}\\
\dss (P\Psi)^{n+1} 
  & = 
    & \dss (P\Psi)^{**} + \frac{\dt}{2} Q\left(\Psi^{n+1}; P^{n+1}\right)
      \label{eq:ImplicitEulerStep}\\
\dss \apsinc P^{n+1} 
  & = 
    & \dss \apsinc P^{n} - \dt\, \nabla\cdot(P\vv)^{n+\half}\,.
    \label{eq:PCompletion}
\end{IEEEeqnarray}
Here we notice that the homogeneous equations \eq{eq:EulerMass} and \eq{eq:EulerPressure}
for $\rho$ and $P$ are not involved in \eq{eq:ExplicitEulerStep} and
\eq{eq:ImplicitEulerStep}. The updates to $\rho^{n+1}$ and $P^{n+1}$ are
entirely determined by the advection step in \eq{eq:AdvectionStep} and
by the completion of the implicit midpoint discretization of the $P$-equation 
in \eq{eq:PCompletion}. 

Therefore, the updated unknowns in the explicit and 
implicit Euler steps \eq{eq:ExplicitEulerStep} and \eq{eq:ImplicitEulerStep} 
are $(\vu, w, \chiprime)$ only. Nevertheless, in order to obtain 
an appropriate approximation of the Exner pressure gradient needed in the 
momentum equation, an auxiliary implicit Euler discretization of the energy 
equation in perturbation form for $\piprime$ from \eq{eq:EulerPiPrime} is 
used in formulating \eq{eq:ImplicitEulerStep}. See 
section~\ref{ssec:SemiImplicitForcing} for details. 

After completion of the steps in \eq{eq:TimeIntegrator} we have 
two redundancies in the thermodynamic variables. In addition to the
primary variables $(\rho, P)$, we also have the perturbation inverse
potential temperature, $\chi'$, and the Exner pressure increment $\pi'$. 
Removal of these redundancies is discussed in 
section~\ref{ssec:Synchronization} below.

Note that the implicit trapezoidal step \eq{eq:TimeIntegrator} and, to a 
lesser extent the treatment of the $P$ in \eq{eq:HalfTimePredictorAdvection}, 
\eq{eq:HalfTimePredictorDivControl}, and \eq{eq:PCompletion}, closely resemble the EULAG/FVM forward-in-time discretization from 
\cite{SmolarkiewiczMargolin1997,PrusaEtAl2008,SmolarkiewiczEtAl2014,SmolarkiewiczEtAl2016, KuehnleinEtAl2019}. 

To avoid misinterpretations, we emphasize that 
\eq{eq:ExplicitEulerStep}-\eq{eq:ImplicitEulerStep} are \emph{not} a
variant of Strang's operator splitting strategy \cite{Strang1968}. To
achieve second-order accuracy, Strang splitting requires all substeps of 
the split algorithm to be second-order accurate individually, aside from
being applied in the typical alternating sequence. This condition is not 
satisfied here as the initial explicit and final implicit Euler steps are 
both only first-order accurate. As shown in \cite{SmolarkiewiczMargolin1993}, 
second-order accuracy results here from a structurally different cancellation 
of truncation errors: By interleaving the Euler steps \eq{eq:ExplicitEulerStep} 
and \eq{eq:ImplicitEulerStep} with one full time step of a second-order 
advection scheme in \eq{eq:AdvectionStep}, one effectively applies the
implicit trapezoidal (or Crank-Nicolson) discretization 
\emph{along the Lagrangian trajectories} described by the advection scheme, 
and this turns out to be second-order accurate, if the trajectories -- the advection step -- are so. 


\section{Discretization details}
\label{sec:DiscretizationDetails}


\subsection{Cartesian grid arrangement}
\label{ssec:GridArrangement}

\begin{figure}
\centering
 \includegraphics[width=.75\columnwidth]{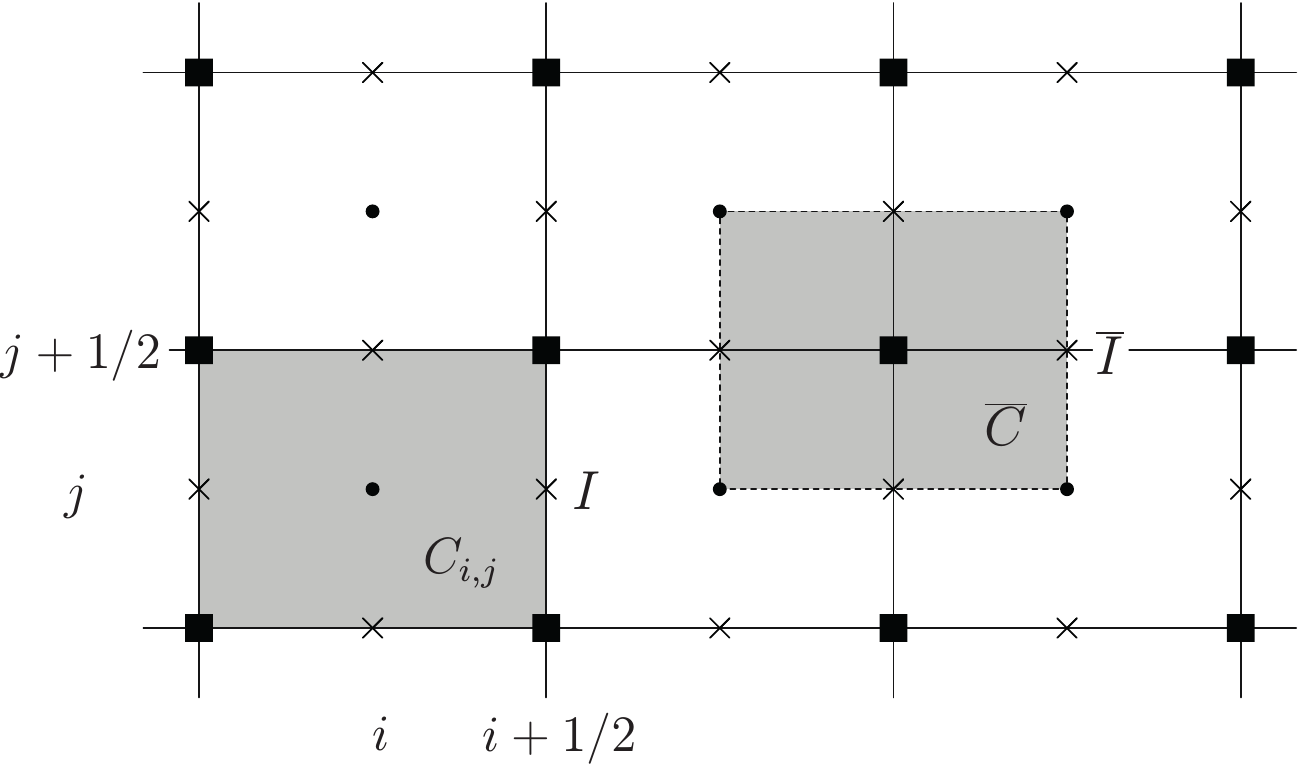}
\caption{Cartesian grid arrangement for two space dimensions. 
$C_{i,j}$: primary finite volumes, 
$\ \bullet\ $: primary cell centers, $I$: primary cell interfaces,
$\times$: centers of both primary and dual cell interfaces, 
$\overline{C}$: dual cells for nodal pressure computation, $\ \rule[1pt]{4pt}{4pt}\ $:
dual cell centers, $\overline{I}$: dual cell interfaces.
\label{fig:GridArrangement}}
\end{figure}

The space discretization of the present scheme for the primary
and auxiliary solution variables 
\begin{equation}
\Sol = \left(\rho, \rho\vu, \rho w, P, P\chiprime \right)^T
\end{equation}
is centered on control volumes $C_{i,j,k}$ formed by a Cartesian mesh with constant,
but not necessarily equal, grid spacings $\dx, \dy, \dz$, and grid indices 
$i = 0, ..., I-1$, $j = 0, ..., J-1$, $k = 0, ..., K-1$ in the three coordinate directions (Figure \ref{fig:GridArrangement} shows a two dimensional $x$-$y$ slice). 
The discrete numerical solution consists of approximate grid cell averages 
\begin{equation}
\Sol_{i,j,k}^n \approx 
\frac{1}{\dx \dy \dz}\int\limits_{C_{i,j,k}} \Sol(\vx,t^n)\, d^3\vx\,.
\end{equation}
The scheme is second-order accurate, so that we can interchangeably interpret 
$\Sol_{i,j,k}^n$ as the cell average or as a point value of $\Sol$ at the center of 
mass of a cell within the approximation order.  

Advection of the specific variables $\Psi$ defined in \eq{eq:PsiDefinition} is mediated 
by staggered-grid components of the advective flux field $(P\vv)^{n+\half}$ referred
to in section~\ref{ssec:TimeDiscretizationOverview} 
above. Specifically, the fluxes
\begin{equation}\label{eq:DiscreteAdvectiveFluxes}
(Pu \Psi)^{n+\half}_{i+\half,j,k}\,,
\qquad
(Pv \Psi)^{n+\half}_{i,j+\half,k}\,, 
\qquad
(Pw \Psi)^{n+\half}_{i,j,k+\half}\,,
\end{equation}
are defined on cell faces $I_{i+\half,j,k}, I_{i,j+\half,k}$, and $I_{i,j,k+\half}$ 
(Figure \ref{fig:GridArrangement}). 
Given the advecting fluxes, e.g., in the $x$-direction 
$(Pu \Psi)^{n+\half}_{i+\half,j,k}$, the associated cell face values 
$\Psi^{n+\half}_{i+\half,j,k}$ 
are determined by a monotone upwind scheme for conservation laws (MUSCL) following 
\cite{vanLeer2006} as described below.


\subsection{Advection}
\label{ssec:Advection}

Any robust numerical scheme capable of performing advection of a scalar in
compressible flows is a good candidate for the generic discrete
advection operators $\mathcal{A}_{\text{1st}}^{\dt}$ and 
$\mathcal{A}_{\text{2nd}}^{\dt}$ introduced in \eq{eq:HalfTimePredictorAdvectionA}
and \eq{eq:AdvectionStep}. 
The present implementation is based on a directionally split 
monotone upwind scheme for conservation laws (MUSCL, see, e.g.,
\cite{vanLeer2006}): 

Suppose the half-time predictor step from \eq{eq:HalfTimePredictorFluxCorrection}, 
the details of which are given in 
section~\ref{sssec:DivControlledAdvectiveFluxes} 
below, has been completed. Then, 
the components of the advecting fluxes $(P\vv)^{n+\half}$ at grid cell faces 
have become available as part of this calculation. Given these fluxes, the advection 
step in \eq{eq:AdvectionStep} is discretized via Strang splitting, so that
\begin{equation}\label{eq:AdvectionStrangSplitting}
\Sol_{i,j}^{**} 
=
\mathcal{A}_{2\text{nd}}^{\dt} \Sol_{i,j,k}^{*} 
\equiv
\mathcal{A}^{\frac{\dt}{2}}_x 
\mathcal{A}^{\frac{\dt}{2}}_y 
\mathcal{A}^{\frac{\dt}{2}}_z 
\mathcal{A}^{\frac{\dt}{2}}_z 
\mathcal{A}^{\frac{\dt}{2}}_y 
\mathcal{A}^{\frac{\dt}{2}}_x\, \Sol_{i,j}^{*} \,,
\end{equation}
where, dropping the indices of the transverse directions 
for simplicity, we have, e.g., 
\begin{equation}
\mathcal{A}^{\frac{\dt}{2}}_x\, \Sol_{i} 
= \Sol_{i}
- \frac{\dt}{2\dx} 
  \left((Pu)^{n+\half}_{i+\half}\, \Psi_{i+\half} 
      - (Pu)^{n+\half}_{i-\half}\, \Psi_{i-\half} \right)
\end{equation}
with 
\begin{IEEEeqnarray}{rCl}\label{eq:AdvSpecifics}
\Psi_{i+\half} 
  & = 
    & \sigma \Psi^{-}_{i+\half} + (1-\sigma)\Psi^{+}_{i+\half}\,,
      \IEEEyesnumber\IEEEyessubnumber*\\[8pt]
\sigma 
  & = 
    & \text{sign}\left((Pu)^{n+\half}_{i+\half}\right)\,,
      \\[8pt]
\Psi^{-}_{i+\half} 
  & = 
    & \Psi_{i} + \frac{\dx}{2} \left(1 - C^{n+\half}_{i+\half} \right)  s_{i}\,,
      \\[8pt]
\Psi^{+}_{i+\half} 
  & = 
    & \Psi_{i+1} - \frac{\dx}{2} \left(1 + C^{n+\half}_{i+\half} \right)  s_{i+1}\,,
      \\[8pt]
C^{n+\half}_{i+\half}
  & =
    & \frac{\dt}{\dx} \frac{(Pu)^{n+\half}_{i+\half}}{(P_{i} + P_{i+1})/2}\,,
      \label{eq:AdvSpecificsCourantNo}\\[8pt]
s_{i}
  & =
    & \text{Lim}
      \left(\frac{\Psi_{i}-\Psi_{i-1}}{\dx}, \frac{\Psi_{i+1}-\Psi_{i}}{\dx}\right)\,,
      \label{eq:AdvSpecificsSlopes}
\end{IEEEeqnarray}
where $P_{i}$ in \eq{eq:AdvSpecificsCourantNo} denotes the fourth component 
of $\Sol_{i}$, and $\text{Lim}(a,b)$ is a slope limiting function 
(see, e.g., \cite{Sweby1984}).

Importantly, the advecting fluxes $(P\vv)^{n+\half}$ are maintained unchanged 
throughout the Strang splitting cycle \eq{eq:AdvectionStrangSplitting}.

The first-order accurate advection operator $\mathcal{A}_{1\text{st}}^{\dt}$
used in \eq{eq:HalfTimePredictorAdvection} is a simplified version
of the above in that the advective fluxes are approximated at the old time level,
i.e., the cell-to-face interpolation formulae for the advective fluxes described
in section~\ref{sssec:DivControlledAdvectiveFluxes} below are evaluated with the
components of $(P\vv)^n$. Optionally, one may also use simple, i.e., not Strang, 
splitting for the advection step of this predictor. In the test shown below, 
we have used the double Strang sweep throughout.


\subsection{Semi-implicit integration of the forcing terms}
\label{ssec:SemiImplicitForcing}

The generalized forcing terms on the right-hand side of \eq{eq:EulerP} are 
discretized in time by the implicit trapezoidal rule. This requires an explicit 
Euler step at the beginning and an implicit Euler step at the end of a time step. 
The implicit Euler scheme is also used to compute the fluxes $(P\vv)^{n+\half}$ 
at the half time level as needed for the advection substep. Below we summarize this 
implicit step in a temporal semi-discretization, explain how this step is used
to access the hydrostatic and pseudo-incompressible balanced models seamlessly, 
provide the node-based spatial discretization, and explain how the 
divergence-controlled momenta are used to generate divergence controlled advective 
fluxes across the faces of the primary control volumes.


\subsubsection{Implicit Euler step and access to hydrostatic and soundproof dynamics}
\label{sssec:ImplicitEuler}

Both $\rho$ and $P$ are frozen in time in this split step because their evolution equations \eq{eq:EulerPMass} and \eq{eq:EulerPP} do not carry a 
right hand side. Hence, the linearized equations including the 
auxiliary potential temperature perturbation equation \eq{eq:EulerPChiPrime}
as well as the hydrostatic and pseudo-incompressible switches, $\ahydro$ and 
$\apsinc$ may be written as 
\begin{IEEEeqnarray}{rCl}\label{eq:LinearizedNonAdvectiveSystem}
U_t
  & = 
    & - c_p (P\Theta)^{\circ} \pi'_x + f V
      \IEEEyesnumber\IEEEyessubnumber*\label{eq:LinearizedNonAdvectiveSystemU}\\[7pt]
V_t
  & = 
    & - c_p (P\Theta)^{\circ} \pi'_y - f U
      \label{eq:LinearizedNonAdvectiveSystemV}\\[0pt]
\ahydro\, W_t
  & =
    & - c_p (P\Theta)^{\circ} \pi'_z + g \frac{\Thetatilde}{\Thetabar}
      \label{eq:LinearizedNonAdvectiveSystemW}\\
\Thetatilde_t
  & =
    & - W\frac{d\Thetabar}{dz}
      \label{eq:LinearizedNonAdvectiveSystemTheta}\\
\apsinc \left(\frac{\partial P}{\partial \pi}\right)^{\circ}\pi'_t
  & =
    & - U_x - V_y - W_z\,,
    \label{eq:LinearizedNonAdvectiveSystemPi}
\end{IEEEeqnarray}
where we have abbreviated
\begin{equation}
(U,V,W,\Thetatilde) = (P u, P v, P w, P(1/\chi)')\,,
\end{equation}
and where the coefficients $(P\Theta)^{\circ}$ and $(\partial P/\partial \pi)^{\circ}$ 
are either those values available when the routine solving the implicit Euler step is called
or they can be adjusted nonlinearly in an outer iteration loop as described in a similar 
context in \cite{SmolarkiewiczEtAl2014}. For all the results shown in this paper we have
used the simpler variant without an outer iteration.

The implicit Euler semi-discretization of \eq{eq:LinearizedNonAdvectiveSystem} 
in time then reads
\begin{IEEEeqnarray}{rCl}\label{eq:LinearizedNonAdvectiveImplicitEuler}
& U^{n+1}
   = 
     U^{n} - \dt \left( c_p (P\Theta)^{\circ} {\pi'}^{n+1}_x - f V^{n+1}\right)
      \IEEEyesnumber\IEEEyessubnumber*\\[7pt]
& V^{n+1}
   = 
     V^n - \dt \left( c_p (P\Theta)^{\circ} {\pi'}^{n+1}_y + f U^{n+1} \right)
      \\[0pt]
& \ahydro\, W^{n+1}
   =
      \ahydro\, W^n - \dt 
       \left( c_p (P\Theta)^{\circ} {\pi'}^{n+1}_z - g \frac{{\Thetatilde}^{n+1}}{\Thetabar} 
       \right)
      \\
& \hspace{-3cm} {\Thetatilde}^{n+1}
   =
     {\Thetatilde}^{n} - \dt\frac{d\Thetabar}{dz}\, W^{n+1}
     & \\
& \hspace{-3cm}\apsinc \left(\frac{\partial P}{\partial \pi}\right)^{\circ} {\pi'}^{n+1}
   =
     \apsinc \left(\frac{\partial P}{\partial \pi}\right)^{\circ} {\pi'}^{n}\label{eq:LinearizedNonAdvectiveImplicitEulerPi}\\ 
   &\qquad\qquad\qquad- \dt \left( U^{n+1}_x + V^{n+1}_y + W^{n+1}_z \right)\,.\nonumber
    \end{IEEEeqnarray}
Straightforward manipulations yield the new time 
level velocity components, 
\begin{IEEEeqnarray}{rCl}\label{eq:LinearizedNonAdvectiveImplicitEulerSolI}
\left(
\begin{array}{c}
U \\ V
\end{array}
\right)^{n+1} 
  & =
    &  \frac{1}{1+(\dt f)^2}
      \left[
      \left(
        \begin{array}{c}
        U + \dt\, f\, V \\ 
        V - \dt\, f\, U
        \end{array}
      \right)^n
\right.      \IEEEyesnumber\IEEEyessubnumber*\\
  & & \left.  - \dt \, c_p (P\Theta)^{\circ} 
      \left(
        \begin{array}{c}
        \pi'_x + \dt\, f\, \pi'_y \\ 
        \pi'_y - \dt\, f\, \pi'_x
        \end{array}
      \right)^{n+1}
      \right] \quad
      \IEEEnonumber\IEEEnosubnumber*\\
W^{n+1} 
  & =
    & \left(\frac{\ahydro W + \dt\, g {\Thetatilde}/\Thetabar}{\ahydro + (\dt\, N)^2}\right)^n \\
  & &     - \dt \, \frac{c_p (P\Theta)^{\circ}}{\ahydro + (\dt\, N)^2} \ {\pi'}^{n+1}_z \,,\IEEEnonumber\IEEEnosubnumber* 
\end{IEEEeqnarray}
with the buoyancy frequency, $N$, given by
\begin{equation}
N^2 = g \frac{1}{\Thetabar}\frac{d\Thetabar}{dz}\,.
\end{equation}
Insertion of the expressions in \eq{eq:LinearizedNonAdvectiveImplicitEulerSolI}
into the pressure equation \eq{eq:LinearizedNonAdvectiveImplicitEulerPi} 
leads to a closed Helmholtz-type equation for ${\pi'}^{n+1}$, 
 \begin{multline}\label{eq:helmholtz}
 \apsinc \left(\dfrac{\partial P}{\partial\pi}\right)^{\circ}\pi'^{n+1} - \dt^2 
\left\{ 
      \left[
        \dfrac{c_p \left(P\Theta\right)^{\circ}}{1+(\dt f)^2} 
          \left(\pi'^{n+1}_x+\dt f\pi'^{n+1}_y
          \right)
        \rule{0pt}{12pt}\right]_x \right.\\
        \left.
       +\left[
          \dfrac{c_p \left(P\Theta\right)^{\circ}}{1+(\dt f)^2} 
            \left(\pi'^{n+1}_y-\dt f\pi'^{n+1}_x
            \right)
         \rule{0pt}{12pt}\right]_y
  + \left[\dfrac{c_p \left(P\Theta\right)^{\circ}}{\ahydro + (\dt N)^2} \ \pi'^{n+1}_z
      \rule{0pt}{12pt}\right]_z
    \right\}
  = R^n
 \end{multline}
with the right-hand side:
\begin{multline}\label{eq:helmholtz_rhs}
R^n= \apsinc \left(\dfrac{\partial P}{\partial\pi}\right)^{\circ}\pi'^n-
\dt
\left\{
\left[\dfrac{U^n +\dt f V^n}{1+(\dt f)^2}\right]_x\right.\\
\left.
+ \left[\dfrac{V^n - \dt f U^n}{1+(\dt f)^2}\right]_y
+\left[\dfrac{\ahydro W^n+\dt\, g\left(\Thetatilde/\Thetabar\right)^n}{\ahydro+(\dt N)^2}\right]_z\right\} \,.
\end{multline}
After its solution, backward re-insertion yields $(U, V, W, \Thetatilde)^{n+1}$.

In all simulations shown in this paper, the Coriolis parameter is set to a
constant, which eliminates the cross-derivative terms $\pi'_{xy}$ from
the elliptic operator in \eq{eq:helmholtz}.

Note that \eq{eq:LinearizedNonAdvectiveImplicitEuler}-\eq{eq:helmholtz_rhs}
reveal how the access to hydrostatic and pseudo-incompressible dynamics is entirely
encoded in the implicit Euler substeps of the scheme, marked by the appearance of 
the $\ahydro$ and $\apsinc$ parameters. In this paper we only demonstrate the behavior 
of the scheme for values of these parameters in $\{0,1\}$, leaving explorations of 
a continuous blending of models with intermediate values of the parameters as well as the development of an analogous switch to geostrophic limiting dynamics to future 
work.


\subsubsection{Pressure gradient and divergence computation in the generalized sources}
\label{sssec:NodalDiscretization}

\begin{figure}
\centering
 \includegraphics[width=.75\columnwidth]{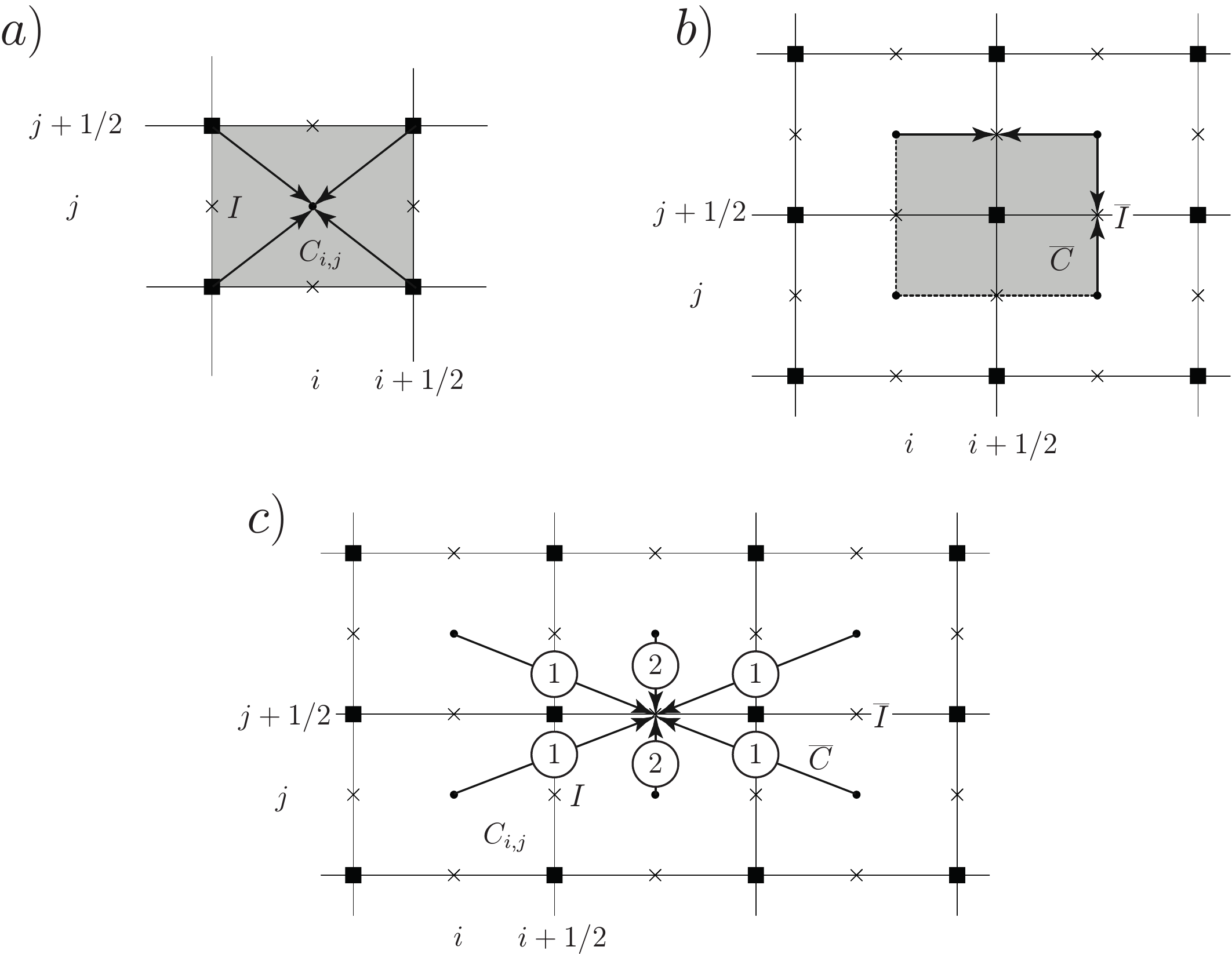}
 \caption{Averaging patterns used in constructing fluxes and cell-centered divergences: 
 a) node-to-cell and analogous cell-to-node averages as in, respectively, \eq{eq:CellToNodeAverage} and 
 \eq{eq:NodeToCellAverage}; 
 b) cell-centered values of flux components $(U,V,W)$ get averaged to the face centers 
 of dual cells in \eq{eq:DualCellFluxAverages}; 
 and c) components 
 of $P\vv$ that are divergence-controlled relative to the nodes are averaged in a 
 particular fashion to cell faces so as to exactly maintain the divergence control. In a) 
 and b) all arrows carry the same weights, so we carry out simple arithmetic averages. In
 c) the numbers in circles indicate relative weights of the participating cell-centered
 values in forming a face value.}
 \label{fig:Interpolations}
\end{figure}

The linearized equations for inclusion of the source terms in
\eq{eq:LinearizedNonAdvectiveSystemU}-\eq{eq:LinearizedNonAdvectiveSystemTheta} 
need to be evaluated at the
cell centers when we apply the two steps of the trapezoidal rule 
in \eq{eq:ExplicitEulerStep} and \eq{eq:ImplicitEulerStep}. To 
this end, the coefficients $(P\Theta)^{\circ}$ are evaluated at the 
cell centers as well, the linearization term from the equation of state
$\left(\partial P/\partial \pi\right)^{\circ}$ is interpolated from 
the cell centers to the nodes according to 
\begin{IEEEeqnarray}{rCl}
a_{i+\half,j+\half,k+\half}
  & =
    & \frac{1}{8}
      \sum_{\lambda,\mu,\nu = 0}^1 a_{i+\lambda,j+\mu,k+\nu} \,,
      \label{eq:CellToNodeAverage}
\end{IEEEeqnarray}
and in a similar way from nodes to cell centers (Figure~\ref{fig:Interpolations}a), 
and the components of the pressure gradient are approximated as
\begin{equation}
\left(\pi'_x\right)_{i,j,k} 
   = 
     \frac{1}{\dx} \left(\piprimehat_{i+\halff,j,k} - \piprimehat_{i-\halff,j,k}\right)
\end{equation}
with
\begin{multline}
\piprimehat_{i+\halff,j,k} = 
    \frac{1}{4}
      \left(\pi'_{i+\halff,j+\halff,k+\halff} + \pi'_{i+\halff,j-\halff,k+\halff}\right.\\
\left.          + \pi'_{i+\halff,j+\halff,k-\halff} + \pi'_{i+\halff,j-\halff,k-\halff}
      \right)\,.
\end{multline}
Analogous formulae hold for the other Cartesian directions. The node-centered 
flux divergence in \eq{eq:LinearizedNonAdvectiveImplicitEulerPi} is formed on the basis 
of the cell-centered components of $\vV = (U,V,W)$, using 
\begin{IEEEeqnarray}{rCl}\label{eq:DualCellFluxAverages}
&\left(U_x\right)_{i+\halff,j+\halff,k+\halff} 
  = 
     \frac{1}{\dx} \left(\Uhat_{i+1,j+\halff,k+\halff} - \Uhat_{i,j+\halff,k+\halff}\right)
      \IEEEyesnumber\IEEEyessubnumber*\\
& \Uhat_{i,j+\halff,k+\halff} 
   = 
     \frac{1}{4}
      \left(U_{i,j+1,k+1} + U_{i,j,k+1} + U_{i,j+1,k} + U_{i,j,k}
      \right)\,,
      \label{eq:Uhat}
\end{IEEEeqnarray}
and analogous formulae for the other Cartesian directions [Figure~\ref{fig:Interpolations}b)].

These spatial discretizations inserted into the temporal semi-discretization
of the implicit Euler step in \eq{eq:LinearizedNonAdvectiveImplicitEuler}
lead to a node-centered discretization of the pressure Helmholtz 
equation based on nine-point and 27-point stencils of the Laplacian in two 
and three dimensions, respectively. The solution provides the required update 
of the node-centered perturbation pressure field and allows us to 
calculate divergence-controlled cell-centered momenta. We note that in the case of the
pseudo-incompressible model, $(\apsinc = 0)$, this amounts to a node-centered 
\emph{exact} projection with a difference approximation that does allow for a 
checkerboard mode in case that the grid has equal spacing in all directions. The authors of \cite{VaterKlein2009} proposed a node-based exact projection that is free of
such modes, but all tests in the present work have used the simpler scheme 
described above.


\subsubsection{Divergence controlled advective fluxes via 
\eq{eq:HalfTimePredictorFluxCorrection}}
\label{sssec:DivControlledAdvectiveFluxes}

Advection is discretized using standard cell-centered flux divergences.  
Thus, the divergence of, e.g., the vector field $\vV = (U,V,W)$ uses the 
discrete approximation
\begin{equation}
\widetilde{\left(U_x\right)}_{i,j,k} 
=
\frac{1}{\dx} \left(U_{i+\halff,j,k} - U_{i-\halff,j,k}\right)\,,
\end{equation}
and analogous expressions for $V_y$ and $W_z$. For stability reasons, we 
need advective fluxes that are divergence-controlled in the sense that they 
are compatible with the Exner pressure 
evolution~\eq{eq:LinearizedNonAdvectiveImplicitEulerPi}. Yet, the Exner pressure 
is stored on grid nodes, so that the flux divergence on the right 
hand side of \eq{eq:LinearizedNonAdvectiveImplicitEulerPi} is node-centered 
but not cell-centered. However, a simple node-to-cell average (Figure \ref{fig:Interpolations}a)
\begin{IEEEeqnarray}{rCl} \label{eq:NodeToCellAverage}
a_{i,j,k}
  & =
    & \frac{1}{8}
      \sum_{\lambda,\mu,\nu = 0}^1 a_{i-\halff+\lambda,j-\halff+\mu,k-\halff+\nu} \,,
\end{IEEEeqnarray}
yields a second-order accurate 
approximation to the cell average. This amounts to approximating the cell-centered 
divergence by the average of the adjacent node-centered divergences. It turns out 
that this is also equivalent to determining the cell-face advective fluxes from the 
interpolation formula
\begin{IEEEeqnarray}{rCl}
U_{i+\halff,j,k} 
  & = 
    & \frac{1}{2} \left(\Uhathat_{i+1,j,k} + \Uhathat_{i,j,k}\right)
      \IEEEyesnumber\IEEEyessubnumber*\\
\Uhathat_{i,j,k} 
  & = 
    & \frac{1}{4}
      \sum_{\mu,\nu = 0}^1 \Uhat_{i,j-\halff+\mu,k-\halff+\nu}
\end{IEEEeqnarray}
with the $\Uhat$ taken from \eq{eq:Uhat}, and with analogous expressions for the other
Cartesian directions. The resulting effective averaging formula takes cell centered
components of $P\vv$ and generates cell face normal transport fluxes (see Figure \ref{fig:Interpolations}c for a two-dimensional depiction).

By this approach, we remove the necessity of separately controlling the advective 
fluxes across the cell faces by a cell-centered elliptic solve (MAC-projection) on the one hand 
and controlling the divergence of the cell-centered velocities by another elliptic 
equation for nodal pressures on the other hand, as in, e.g., \cite{AlmgrenEtAl2006,BellEtAl1989,BenacchioEtAl2014,SchneiderEtAl1999}.
Thus, the present scheme works with the node-based discretization of the Helmholtz 
equation only. We note in passing that this approach requires an \emph{exact} projection
for the nodal divergence. 


\subsection{Synchronization of auxiliary variables}
\label{ssec:Synchronization}

The proposed scheme achieves large time step capabilities, i.a., by introducing
two additional auxiliary variables that are to be synchronized with the current
state represented by the primary cell averages of $(\rho, \rho\vv, P)$ after 
each time step. 


\subsubsection{Adjustment of the potential temperature perturbation}
\label{ssec:PotTempAdj}

This synchronization is straightforward for the inverse of the potential temperature 
$\chi = \chiprime + \chibar$. After completion of the $n$th time step we let
\begin{equation}
\chiprime_{i,j,k}^{n} = \chi_{i,j,k}^{n} - \chibar_{k}^{n}\,,
\end{equation}   
where we have assumed gravity to be aligned with the $z$-coordinate
direction so that the discrete version of $\chibar(z)$ depends on the 
associated index $k$ only. Also, in all simulations in this paper we have 
set $\chibar_k^{n} \equiv \chibar_k^{0}$, i.e., we have not re-computed
the horizontal average $\chibar(z)$ during the simulations.
An alternative option better suited for large-scale long-time simulations 
is to invoke a horizontal, possibly local, averaging procedure to extract 
$\chibar$ from $\chi$ at least every few time steps. We leave testing 
this option to future work.	 


\subsubsection{Synchronization of nodal and cell pressures}
\label{ssec:SynchronizationChi}

In section~\ref{sssec:DivControlledAdvectiveFluxes} we constructed the cell-centered advective flux divergence from the arithmetic
average of the divergences obtained on the adjacent nodes. By the same 
reasoning the cell-centered update of $P$ that results
from these cell-centered divergences corresponds to the node-to-cell average 
\eq{eq:NodeToCellAverage} for 
$(\partial P/\partial \pi)^{\circ} (\pi^{n+1}-\pi^{n})$.
If, in addition, the pressure Helmholtz equation from \eq{eq:helmholtz} is
solved with an outer iteration such that after convergence this coefficient
is approximated by  
\begin{equation}
\left(\frac{\partial P}{\partial \pi}\right)^{\circ}_{i+\halff,j+\halff,k+\halff}
= \left(\frac{P^{n+1} - P^{n}}{\pi^{n+1} - \pi^{n}}\right)_{i+\halff,j+\halff,k+\halff}\,,
\end{equation}
then the cell-centered time updates of $P$ are guaranteed to always equal the
node-to-cell average of their nodal counterparts. As a consequence, a potential
cumulative desynchronization over many time steps of the nodal Exner pressure 
values and the cell-centered values of $P$ is avoided. 

For the tests shown in this paper, we have not used such an outer iteration, yet we did not observe a desynchronization even over tens of thousands of time steps.


\section{Numerical Results}
\label{sec:Results}

The algorithm described in the previous sections was tested on a suite of benchmarks of dry compressible dynamics on a vertical $x-z$ slice at various scales. The suite draws on the set of \cite{Benacchio2014,BenacchioEtAl2014} including a cold air bubble and nonhydrostatic inertia-gravity waves, and adds to it three larger scale configurations for the inertia-gravity waves, with the aim to validate both the robustness and accuracy of the new buoyancy-implicit strategy, and the scheme's capability of accessing compressible, pseudo-incompressible, and hydrostatic dynamics. We remark that the present paper does not focus on efficiency. While the coding framework is 3D-ready, we leave parallelization and performance on three-dimensional tests for future work. The scheme is implemented in plain C language and uses the Bi-CGSTAB linear solver \cite{Vandervorst1992} for the solution of the elliptic problems. The solver tolerance was set at $10^{-8}$ throughout. We also define the advective Courant number as:
\begin{equation}
\textrm{CFL}_\textrm{adv} = \max\limits_{i \in \{1,2,3\}}\left(\frac{\dt v_i}{\dx_i}\right)
\end{equation}
where $v_i$ are the components of the velocity, $\dx_i$ the grid spacing in the $i$ direction, and the acoustic Courant number as:
\begin{equation}
\textrm{CFL}_\textrm{ac} = \max\limits_{i \in \{1,2,3\}}\left[\frac{\dt (v_i+c)}{\dx_i}\right]
\end{equation}
where $c=\sqrt{\gamma RT}$ denotes the speed of sound.

\subsection{Density current}

The first test case, proposed in \cite{StrakaEtAl1993}, concerns the simulation of a falling bubble of cold air in a neutrally stratified atmosphere $(x,z)\in[-25.6,25.6]\times[0,6.4]\,\textrm{km}^2$. The reference potential temperature and pressure  are $\theta_{ref}=300\,\textrm{K}$ and $p_{ref}=10^5\,\textrm{Pa}$, the thermal perturbation is:
\begin{equation}
 T'=\begin{cases}
           0~\textrm{K} \qquad&\textrm{if}\;r>1\\
	   -15\left[1+\cos(\pi r)\right]/2 \;\;\textrm{K}\qquad&\textrm{if}\;r<1   
          \end{cases},
\end{equation} 
where $r=\left\{[(x-x_c)/x_r]^2+[(z-z_c)/z_r]^2\right\}^{0.5}$, $x_c=0~\textrm{km}$, $x_r=4~\textrm{km}$, $z_c=3~\textrm{km}$ and $z_r=2~\textrm{km}$. 
Boundary conditions are solid walls on top and bottom boundaries and periodic elsewhere. In order to obtain a converged solution, artificial diffusion terms $\rho\mu\nabla^2\mathbf{u}$ and
$\rho\mu\nabla^2\Theta$ are added to the momentum and $P$-equations, respectively, with $\mu=75\,\textrm{m$^2$s$^{-1}$}$. The terms are non-stiff, discretized by the explicit Euler method individually, and tied into the scheme via Operator splitting just before the second backward
Euler step \eq{eq:ImplicitEulerStep}.

\begin{figure}
\centering
 \includegraphics[width=.75\columnwidth]{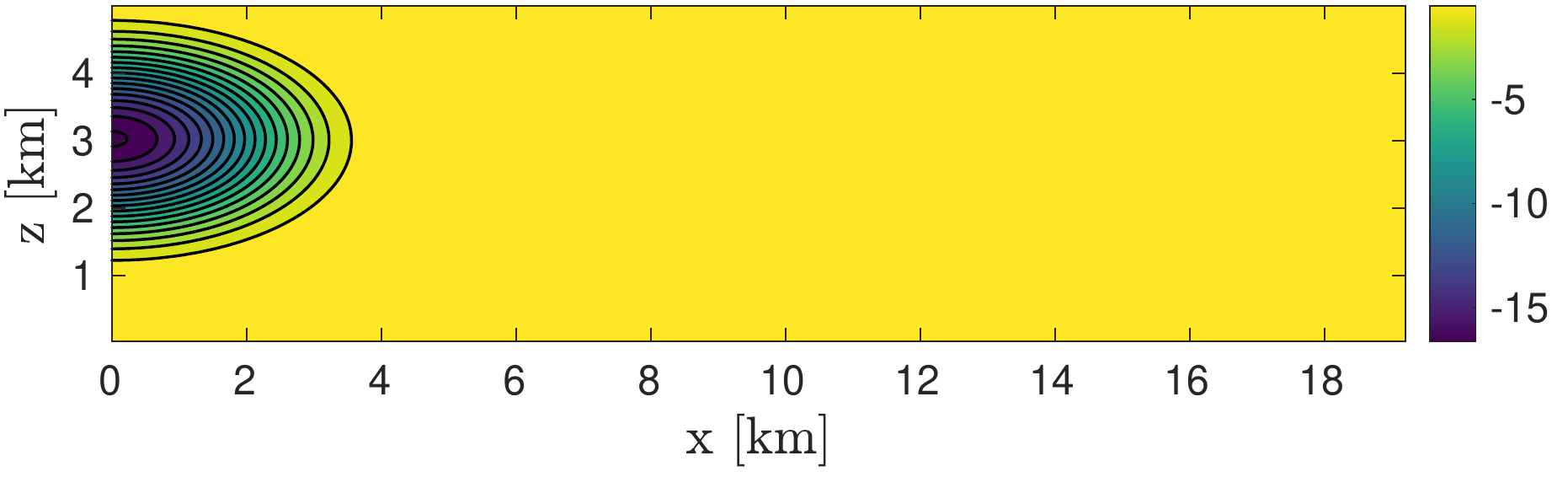}
 \includegraphics[width=.75\columnwidth]{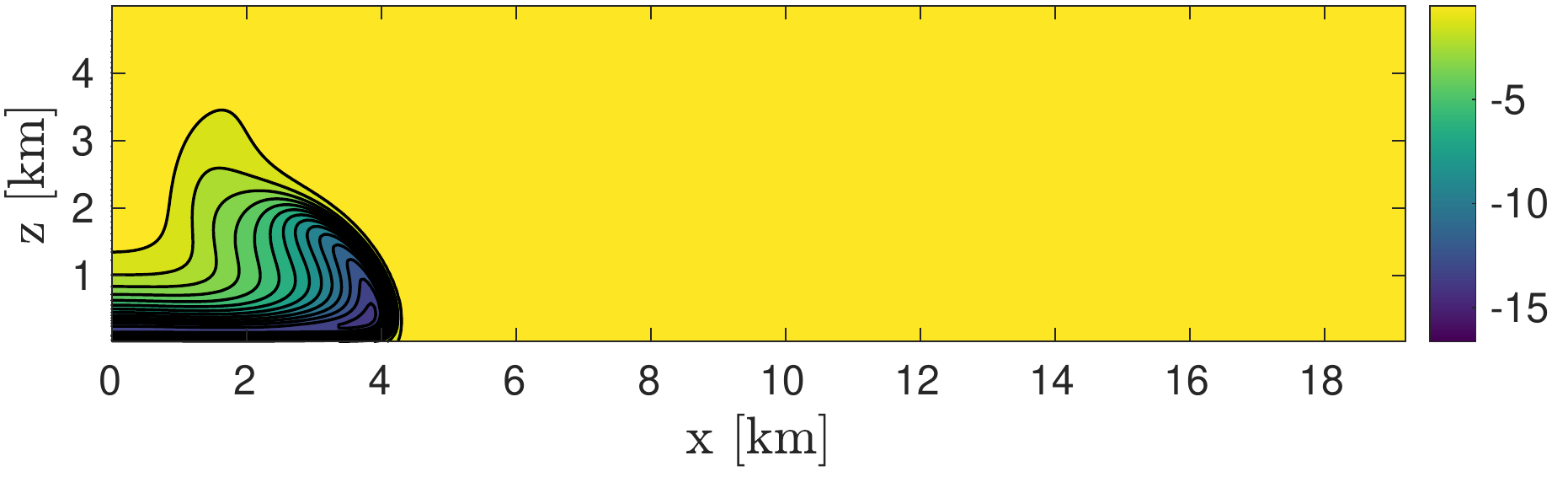}
 \includegraphics[width=.75\columnwidth]{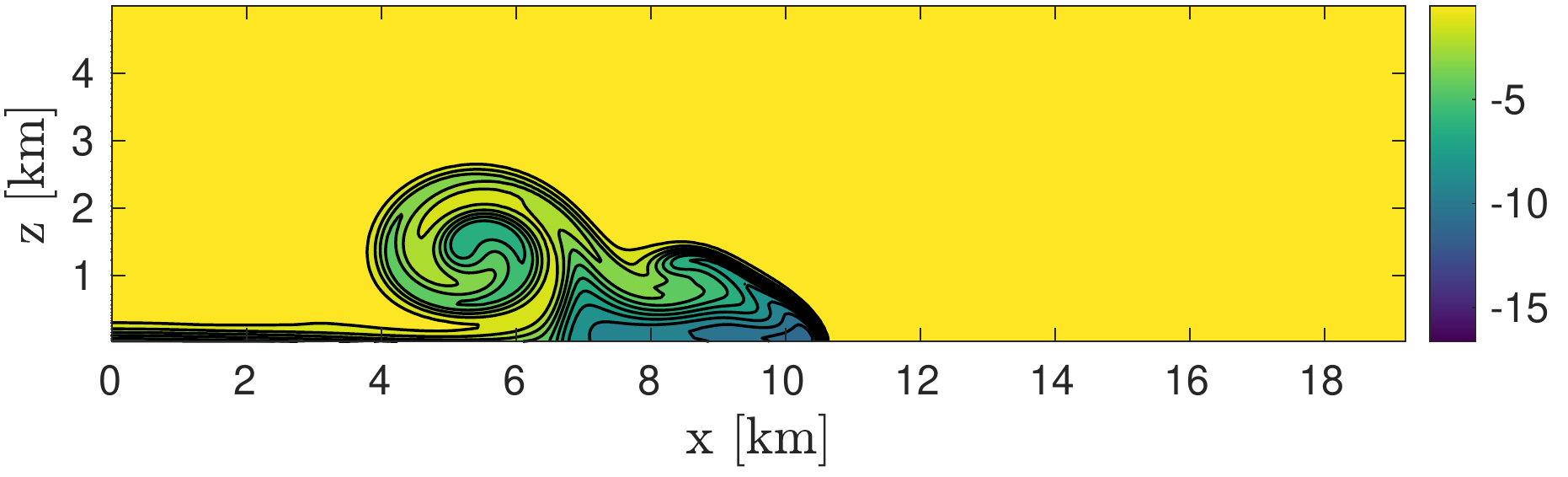}
 \includegraphics[width=.75\columnwidth]{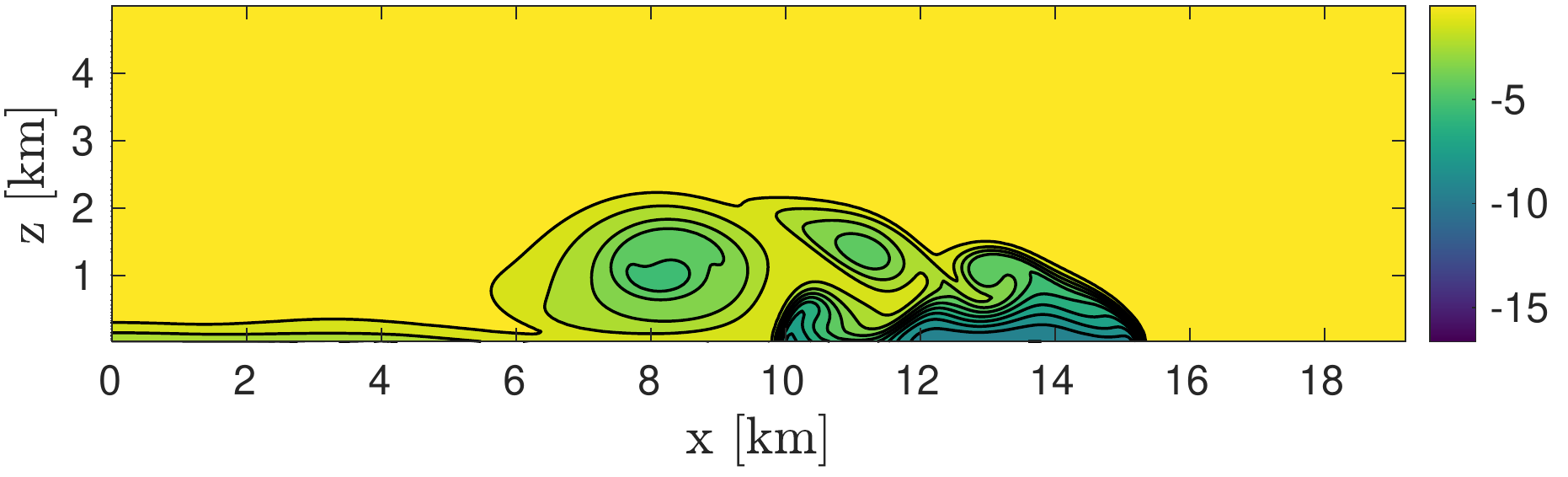}
 \caption{Potential temperature perturbation at times (top to bottom) $t=0,\,300,\,600,\,900\,\textrm{s}$ for the density current test case at spatial resolution $\Delta x=\Delta z=50\,\textrm{m}$, CFL$_\textrm{adv}=0.96$. Contours in the range $[-16.5,\,-0.5]\,\textrm{K}$ with a $1\,\textrm{K}$ contour interval.}
 \label{fig:straka}
\end{figure}

\begin{figure}
\centering
  \includegraphics[width=.75\columnwidth]{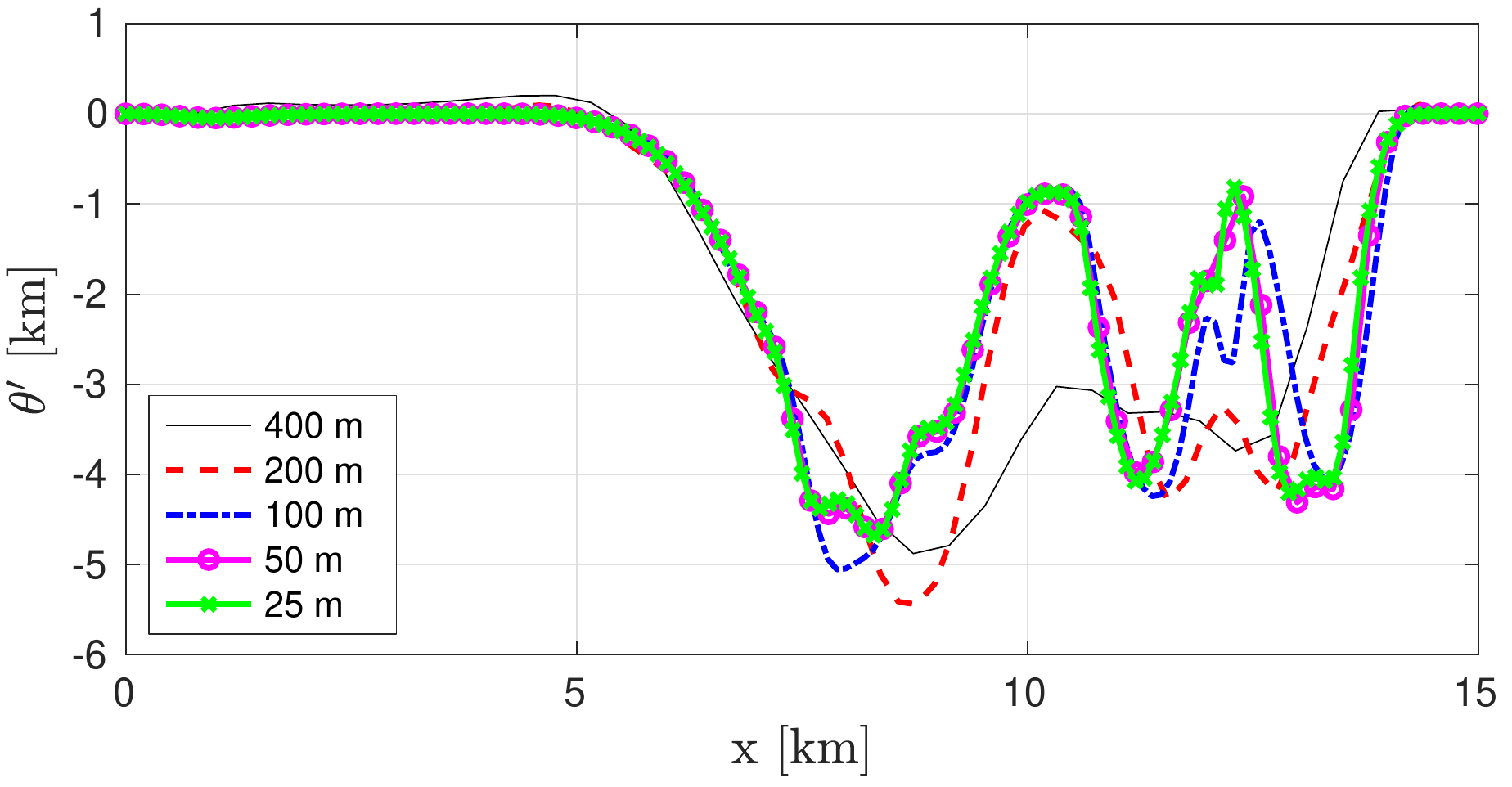}
 \caption{One-dimensional cut at height $z=1200\,\textrm{m}$ for the potential temperature perturbation at final time $t=900\, s$ in the density current test case run with CFL$_\textrm{adv}=0.96$. Spatial resolutions $\Delta x=\Delta z=400\,\textrm{m (thin solid black line)}$, $200\,\textrm{m (dashed red line)}$, $100\,\textrm{m (dashed-dotted blue line)}$, $50\,\textrm{m (thin solid black line)}$, $25\,\textrm{m (thin solid black line)}$.}
 \label{fig:straka1dcut}
\end{figure}

In the reference setup for this case, the buoyancy-implicit model is run at a resolution $\Delta x=\Delta z=50\,\textrm{m}$ with time step chosen according to CFL$_\textrm{adv}=0.96$. Driven by its negative buoyancy, the initial perturbation moves downwards, impacts the bottom boundary and travels sideways developing vortices (Figure \ref{fig:straka}). The numerical solution converges with increasing spatial resolution (Figure \ref{fig:straka1dcut}), and the final perturbation amplitude and front position agree with published results (Table \ref{tab:straka_minmax}, for comparison see, e.g., \cite{GiraldoRestelli2008} and the similar table in \cite{MelvinEtAl2018}). The final minimum potential temperature perturbation at $25\,\textrm{m}$ resolution agrees with the result in \cite{MelvinEtAl2018} up to the third decimal digit.

\begin{table}
\small
\begin{centering}
\begin{tabular}{lccc}
\toprule 
\shortstack{Grid size \\ $[$m$]$}& \shortstack{$\theta'_\textrm{min}$ [K]} & \shortstack{$\theta'_\textrm{max}$ [K]} & \shortstack{Front\\ location [m]} \tabularnewline
\midrule 
400 &  $-8.1466$ & $ 0.2685$ & 14125\tabularnewline
200 &  $-8.9358$ & $0.2294$ & 14884\tabularnewline
100 &  $-9.2154$ & $0.1787$ & 15199\tabularnewline
50 &  $-9.5061$ & $0.0903$ & 15326\tabularnewline
25 &  $-9.6555$ & $0.0138$ & 15381\tabularnewline
\bottomrule 
\end{tabular}

\par\end{centering}
\caption{Minimum and maximum potential temperature perturbation and front location
(rightmost intersection of $-1\,\textrm{K}$ contour with $z=0$) for the density current test at several resolution values.}%
\label{tab:straka_minmax}
\end{table}

\subsection{Inertia-gravity waves} 

The next set of tests consists of gravity waves in a stably stratified channel with constant buoyancy frequency $N=0.01\textrm{s}^{-1}$, $\theta(z=0)=300\,\textrm{K}$, horizontal extension $x\in[0,x_N]$, and vertical extension $z=10\,\textrm{km}$, proposed in \cite{SkamarockKlemp1994}. The thermal perturbation is:
\begin{equation}
 \theta'(x, z, 0)=0.01~\textrm{K}*\dfrac{\sin(\pi z/H)}{1+[(x-x_c)/a]^2}\label{eq: init_theta_pert_igw} 
\end{equation} 
with $H=10~\textrm{km}$, $x_c=100~\textrm{km}$, $a=x_N/60$, and there is a background horizontal flow $u=20~\textrm{m s}^{-1}$. We consider three configurations for the horizontal extension $x_N=300\,\textrm{km},\,6\,000\,\textrm{km},\,48\,000\,\textrm{km}$, with respective final times $T=3\,000\,\textrm{s},\,60\,000\,\textrm{s},\,480\,000\,\textrm{s}$. The first two configurations correspond to the nonhydrostatic case and the hydrostatic case of \cite{SkamarockKlemp1994}, the third planetary-scale configuration is introduced in this paper. In all configurations, the buoyancy-implicit model is run with $300\times10$ cells, as in \cite{SkamarockKlemp1994}, and CFL$_\textrm{adv}=0.9$.

\begin{figure}
\centering
 \includegraphics[width=.75\columnwidth]{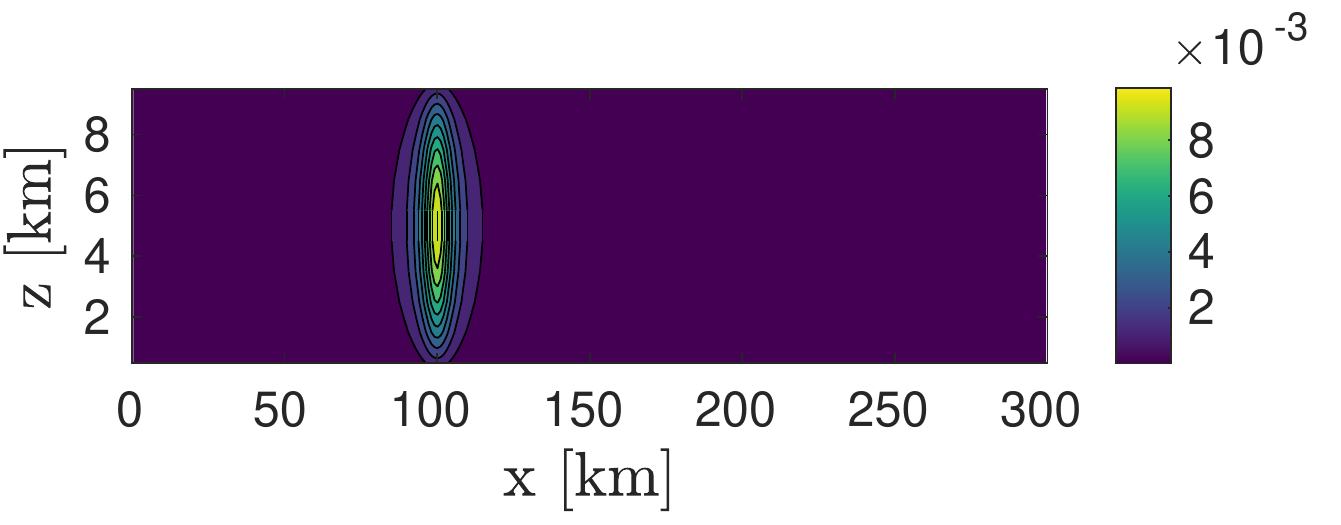}
 \includegraphics[width=.75\columnwidth]{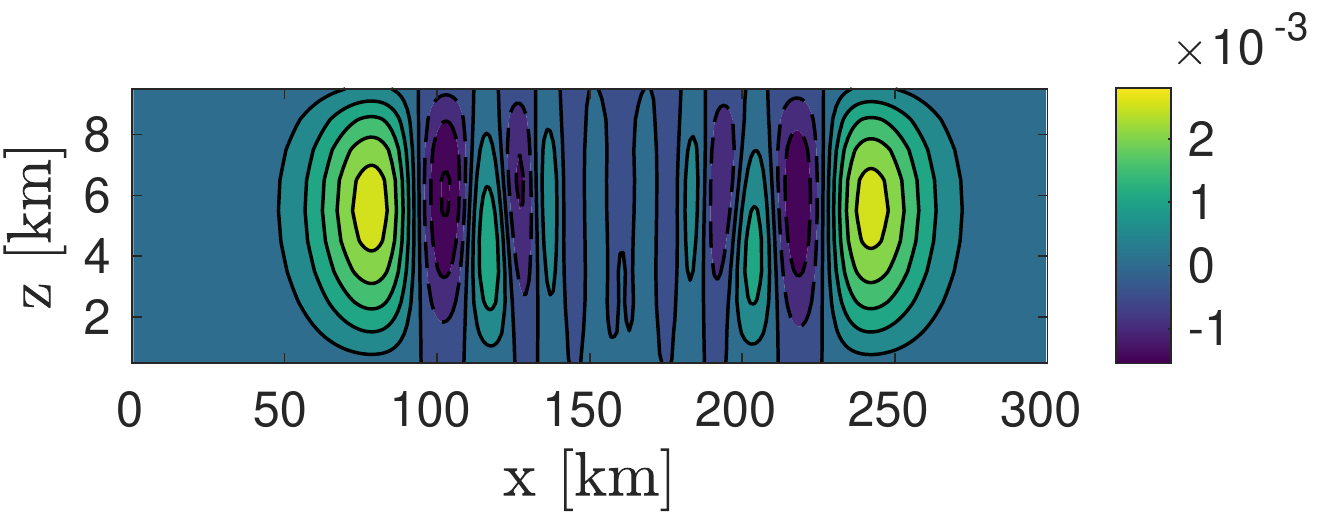}
 \caption{Initial (top) and final (at $T=3000\,\textrm{s}$, bottom) potential temperature perturbation for the nonhydrostatic inertia-gravity wave test from \cite{SkamarockKlemp1994}, $\Delta x=\Delta z=1\,\textrm{km}$, CFL$_\textrm{adv}=0.9$. Contours in the range $[0, 0.01]\,\textrm{K}$ with a $0.001\,\textrm{K}$ interval (top), $[-0.0025, 0.0025]\,\textrm{K}$ with a $0.0005\,\textrm{K}$ interval (bottom). Negative contours are dashed.}
 \label{fig:SK94_NH}
\end{figure}

\begin{figure}
\centering
 \includegraphics[width=.49\columnwidth]{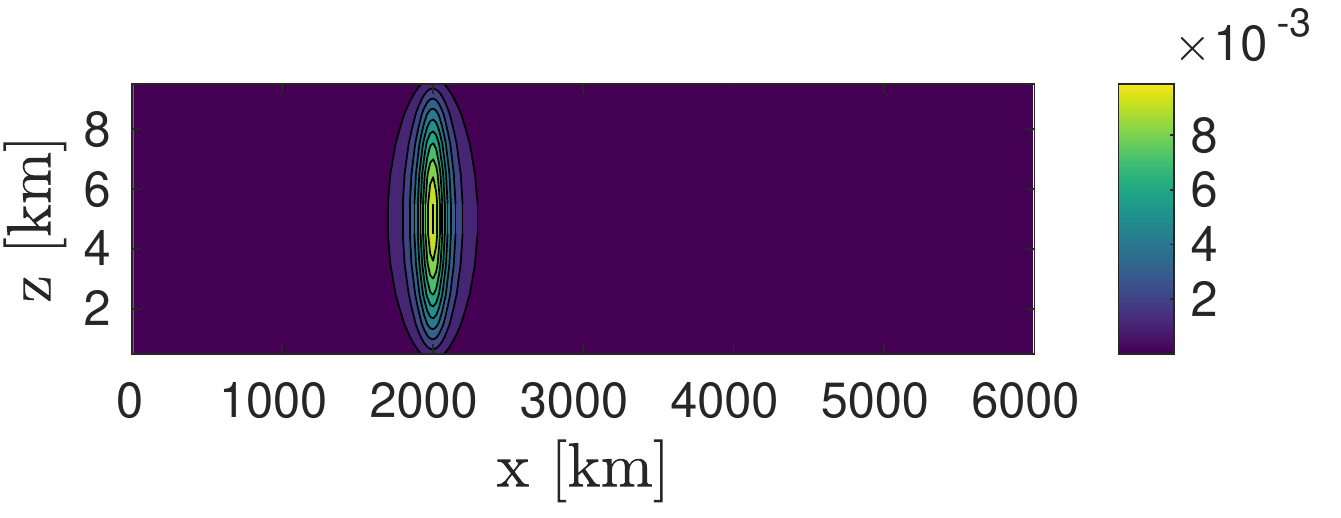}
 \includegraphics[width=.49\columnwidth]{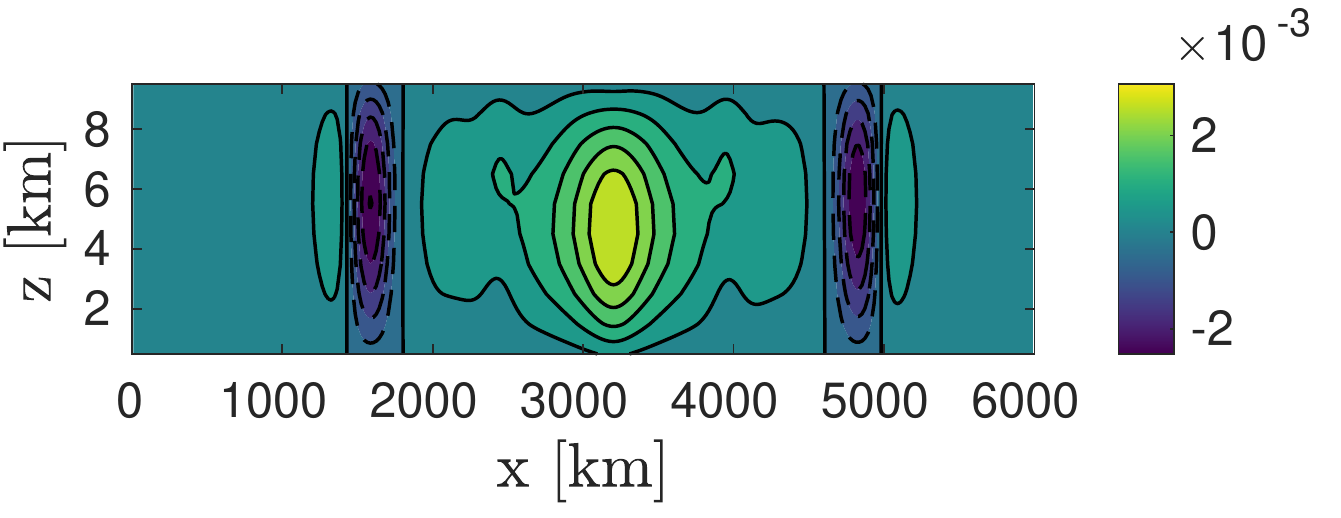}\\
 \includegraphics[width=.49\columnwidth]{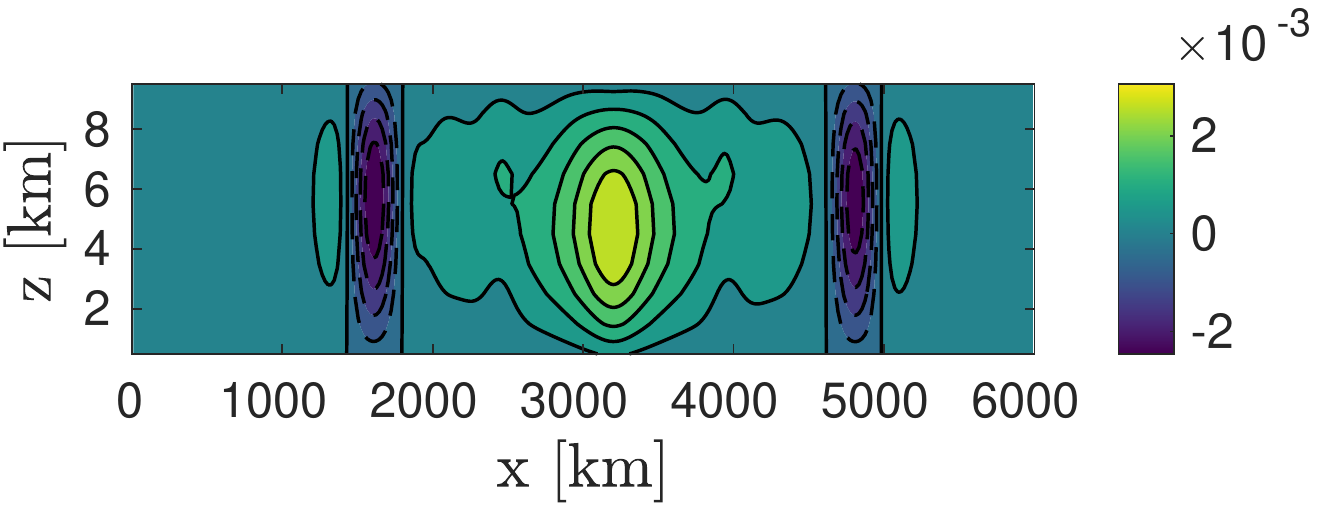}
 \includegraphics[width=.49\columnwidth]{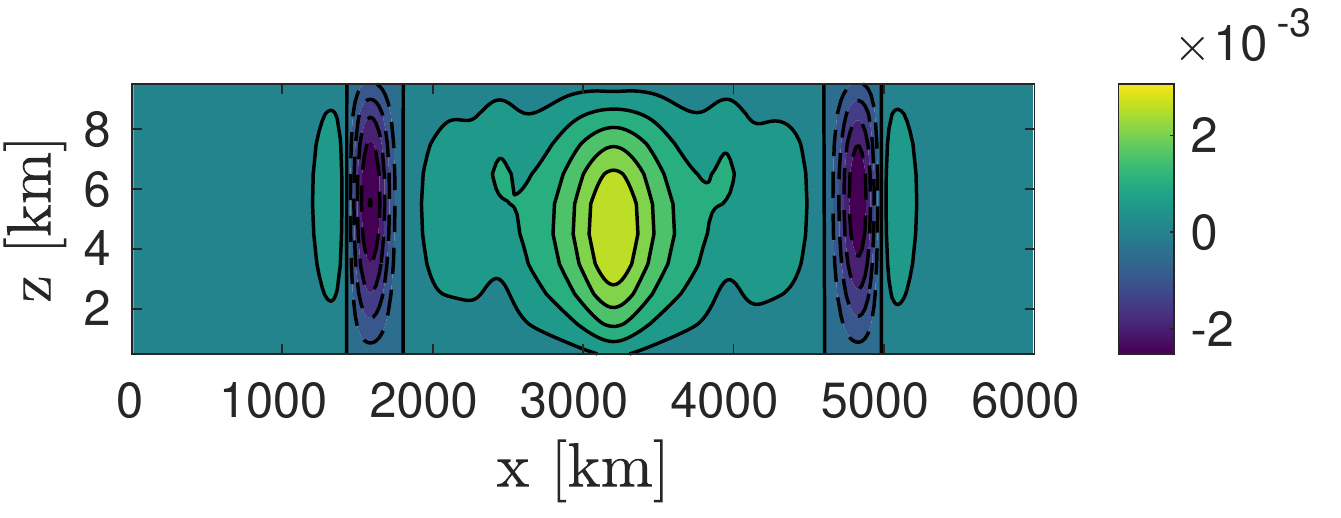}\\
  \includegraphics[width=.49\columnwidth]{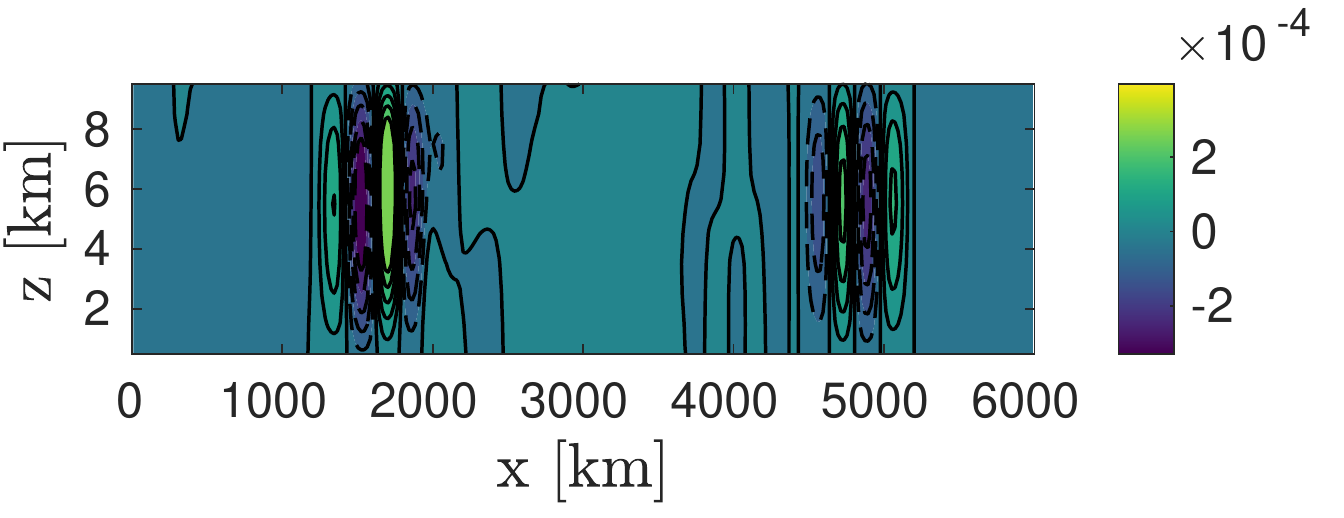}
 \includegraphics[width=.49\columnwidth]{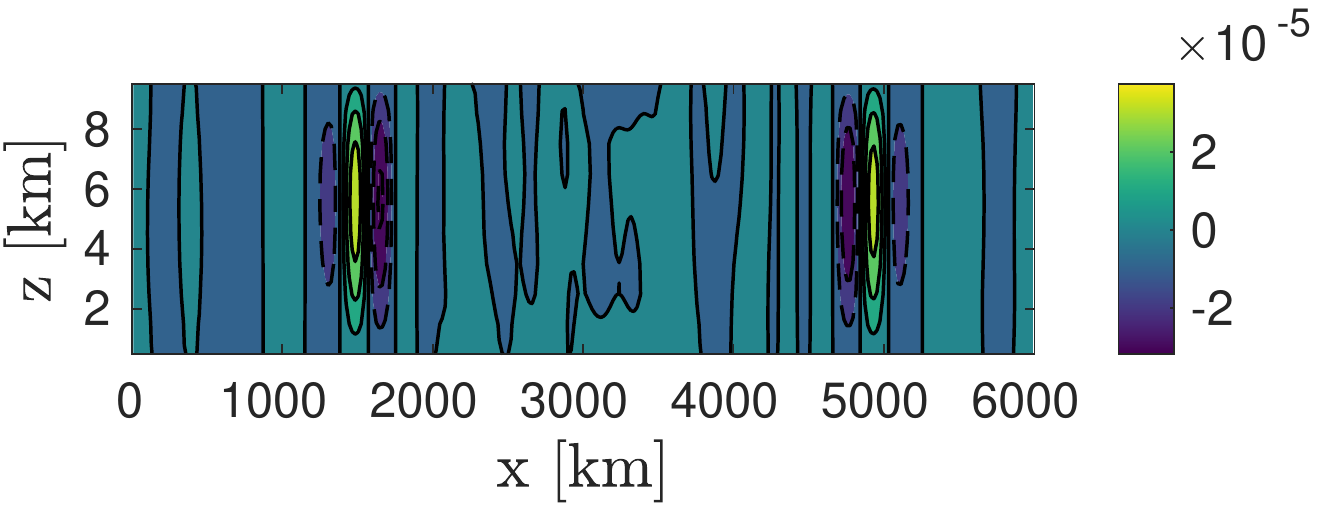}
 \caption{Potential temperature perturbation for the hydrostatic inertia-gravity wave test from \cite{SkamarockKlemp1994},  $\Delta x=20\,\textrm{km},\,\Delta z=1\,\textrm{km}$, CFL$_\textrm{adv}=0.9$. Initial data (top left) and computed value at final time $T=60000\,\textrm{s}$ in compressible mode (top right), pseudo-incompressible mode (middle left), hydrostatic mode (middle right). Contours as in Figure \ref{fig:SK94_NH}. The bottom plots show the difference between the compressible run and the pseudo-incompressible run (left) and between the compressible run and the hydrostatic run (right). Contours in the range $[-2.5, 2.5]*10^{-4}\,\textrm{K}$ with a $5*10^{-5}\,\textrm{K}$ interval (left), $[-5, 5]*10^{-5}\,\textrm{K}$  with a $10^{-5}\,\textrm{K}$ interval (right). Negative contours are dashed.}
 \label{fig:SK94_H} 
 \end{figure}

In the first configuration, the initial perturbation spreads out onto gravity waves driven by the underlying buoyancy stratification (Figure \ref{fig:SK94_NH}). In the second configuration, which is run with rotation (Coriolis parameter value $f=10^{-4}\,\textrm{s}^{-1}$) , a geostrophic mode is also present in the center of the domain (Figure \ref{fig:SK94_H}). In both cases, the values obtained by running the compressible model (COMP) closely resemble published results in the literature including, for the nonhydrostatic case, the buoyancy-explicit compressible result in \cite{BenacchioEtAl2014}. At CFL$_\textrm{adv}=0.9$, the time step used in the first configuration is $\Delta t\approx44.83\,\textrm{s}$, a 12 times larger value than \cite{BenacchioEtAl2014}'s $3.75\,\textrm{s}$. The time step value used here is also in line with \cite{MelvinEtAl2018}, who ran the configuration with $\dt=12\,\textrm{s}$ at buoyancy-implicit CFL$=0.3$. For the second configuration at CFL$_\textrm{adv}=0.9$, the time step used is $\Delta t\approx896.48\,\textrm{s}$, equivalent to an acoustic $\mathrm{CFL}_{\mathrm{ac}}\approx309.5$ and $N\dt = 8.96$.

\begin{figure}
\centering
 \includegraphics[width=.49\columnwidth]{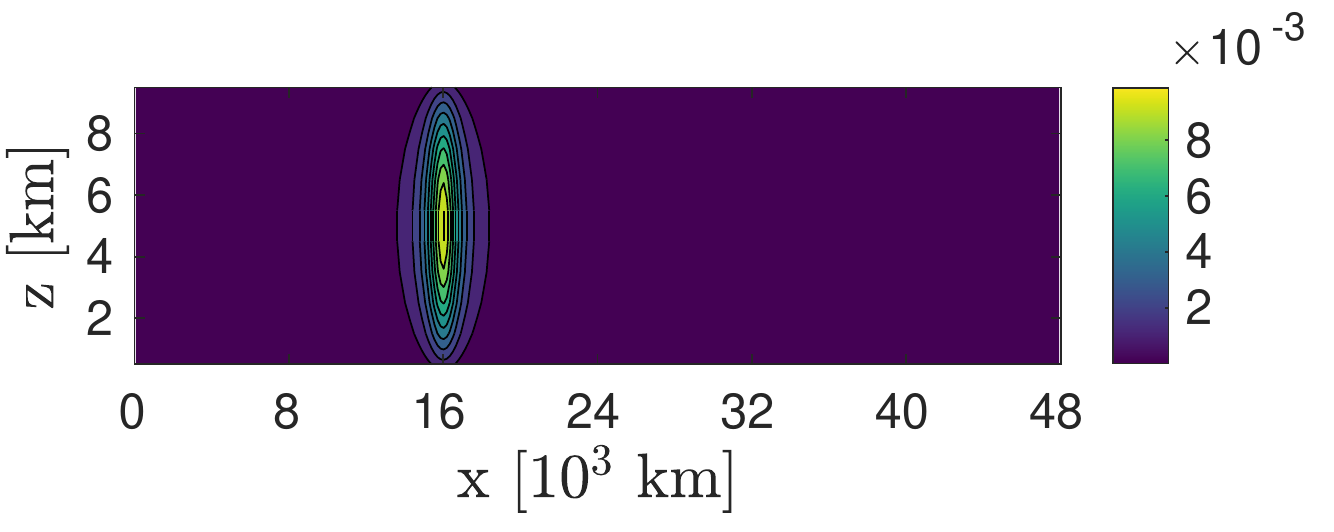}
 \includegraphics[width=.49\columnwidth]{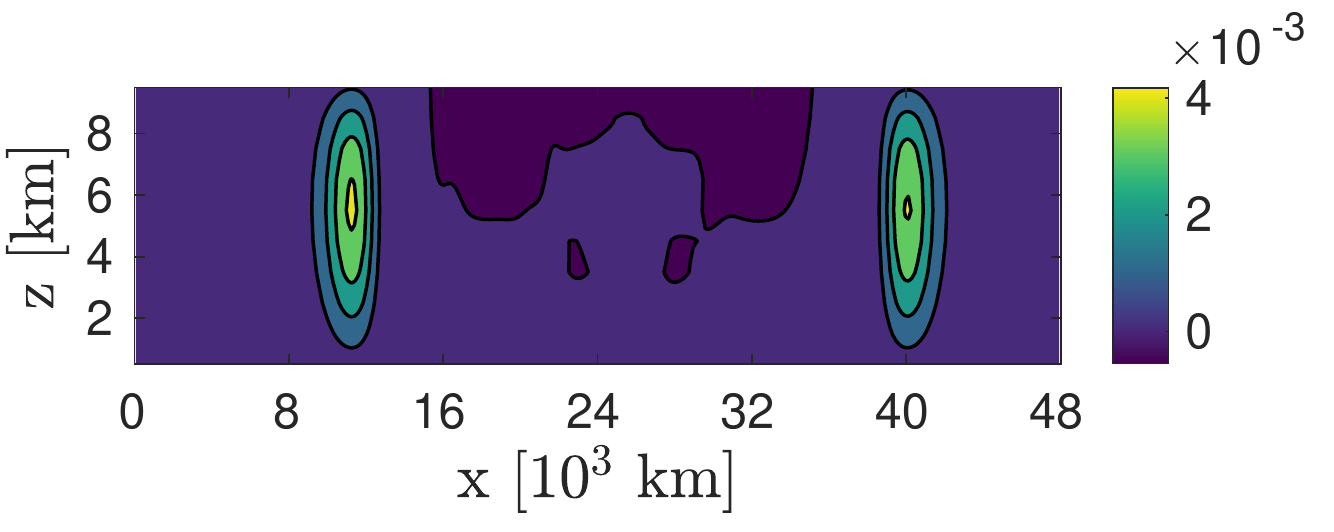}\\
 \includegraphics[width=.49\columnwidth]{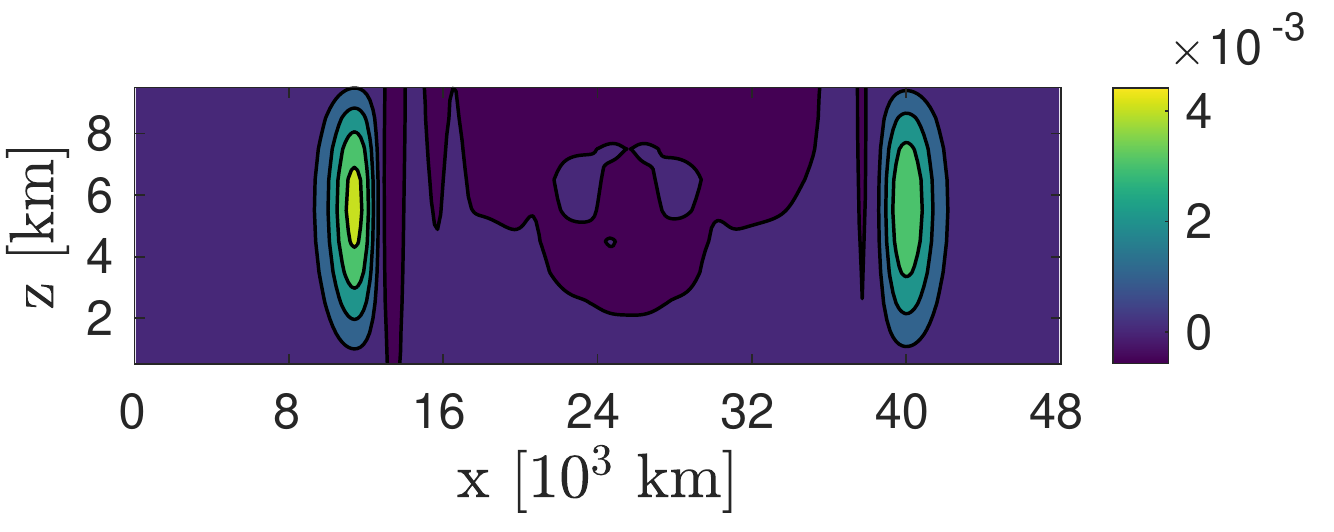}
 \includegraphics[width=.49\columnwidth]{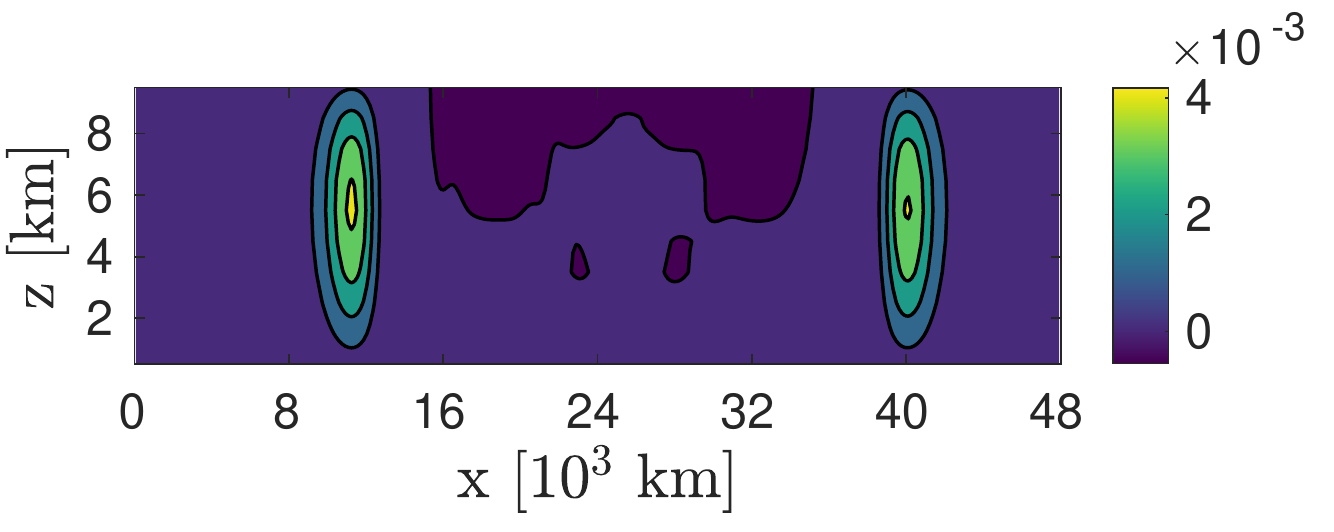}\\
  \includegraphics[width=.49\columnwidth]{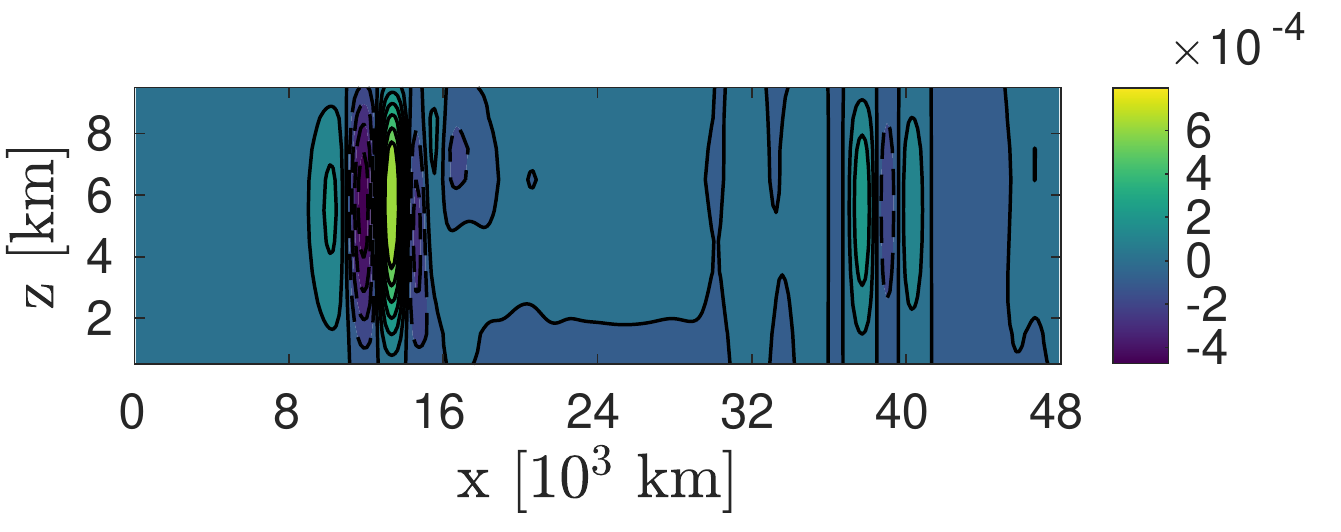}
 \includegraphics[width=.49\columnwidth]{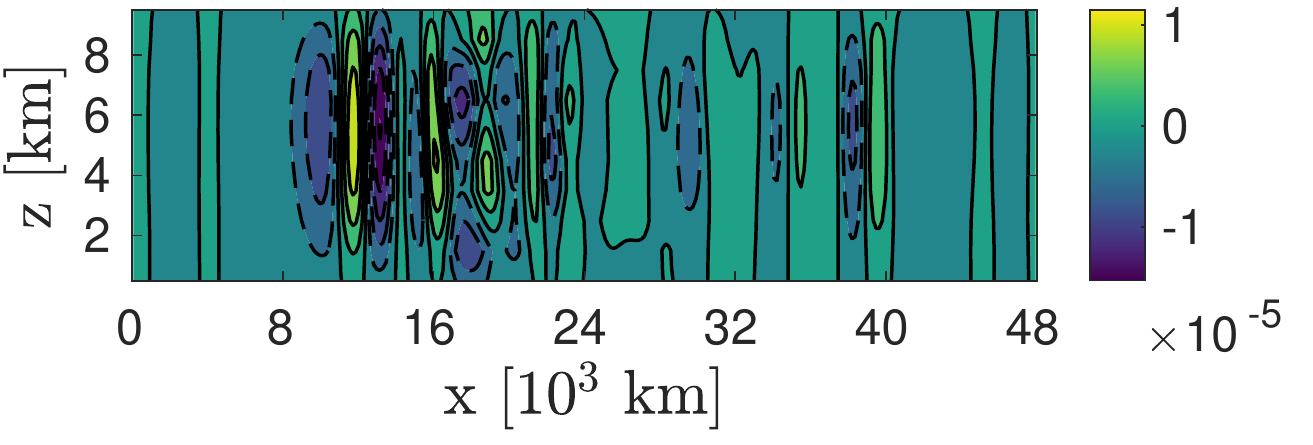}
 \caption{Potential temperature perturbation for the planetary-scale inertia-gravity wave test,  $\Delta x=160\,\textrm{km},\,\Delta z=1\,\textrm{km}$, CFL$_\textrm{adv}=0.9$. Initial data (top left, contours as in Figures \ref{fig:SK94_NH}-\ref{fig:SK94_H}) and computed value at final time $T=480000\,\textrm{s}$ in compressible mode (top right), pseudo-incompressible mode (middle left), hydrostatic mode (middle right). Contours in the range $[-0.005, 0.005]\,\textrm{K}$ with a $0.001\,\textrm{K}$ interval. The bottom plots show the difference between the compressible run and the pseudo-incompressible run (left) and between the compressible run and the hydrostatic run (right). Contours in the range $[-4, 6]*10^{-4}\,\textrm{K}$ with a $10^{-4}\,\textrm{K}$ interval (left), $[-1.5, 1.5]*10^{-5}\,\textrm{K}$ with a $3*10^{-6}\,\textrm{K}$ interval (right). Negative contours are dashed.}
 \label{fig:SK94_P} 
 \end{figure}

The third new planetary-scale configuration is run without rotation to suppress the otherwise dominant geostrophic mode and highlight the wave dynamics. At final time $T=480\,000\,\textrm{s}$ ($\approx5.5$ days), the solution quality with the compressible model is good in terms of symmetry, absence of oscillations, and final position of the outermost crests (Figure \ref{fig:SK94_P}). The time step in this run at $\mathrm{CFL}_\mathrm{adv} = 0.9$, is $\dt\approx7100\,s$, equivalent to $N \dt\approx71$ and to an acoustic $\mathrm{CFL}_\mathrm{ac}\approx2.4\cdot 10^3$.

For the two largest configurations, we also report the pseudo-incompressible (PI) result obtained using $\alpha_P=0$, i.e. by switching off compressibility zeroing the diagonal term in the Helmholtz equation, and the hydrostatic (HY) result obtained using $\alpha_w=0$, i.e. by zeroing the dynamic tendency of the velocity in the vertical momentum equation (middle panels of Figures \ref{fig:SK94_H}-\ref{fig:SK94_P}), together with the differences with the compressible result, COMP$-$PI and COMP$-$HY (bottom panels of Figures \ref{fig:SK94_H}-\ref{fig:SK94_P}). The discrepancies with the compressible result are larger with the pseudo-incompressible model than with the hydrostatic model. Moreover, COMP$-$PI grows with larger horizontal scales and COMP$-$HY shrinks as expected with smaller vertical-to horizontal domain size aspect ratios.

\subsection{Superposition of acoustic-gravity waves and inertia-gravity waves}

As final corroboration of the properties of the model, the hydrostatic configuration is rerun with a different value of the Coriolis parameter $f=1.03126*10^{-4}\,\textrm{s}^{-1}$, initial temperature $T(z=0)=250\,\textrm{K}$, isothermal background distribution, and no background flow. A time step of $\dt=0.125\,\textrm{s}$ is used as in \cite{BaldaufBrdar2013} for a run with $1200\times80$ cells. 

\begin{figure}
\centering
 \includegraphics[width=.75\columnwidth]{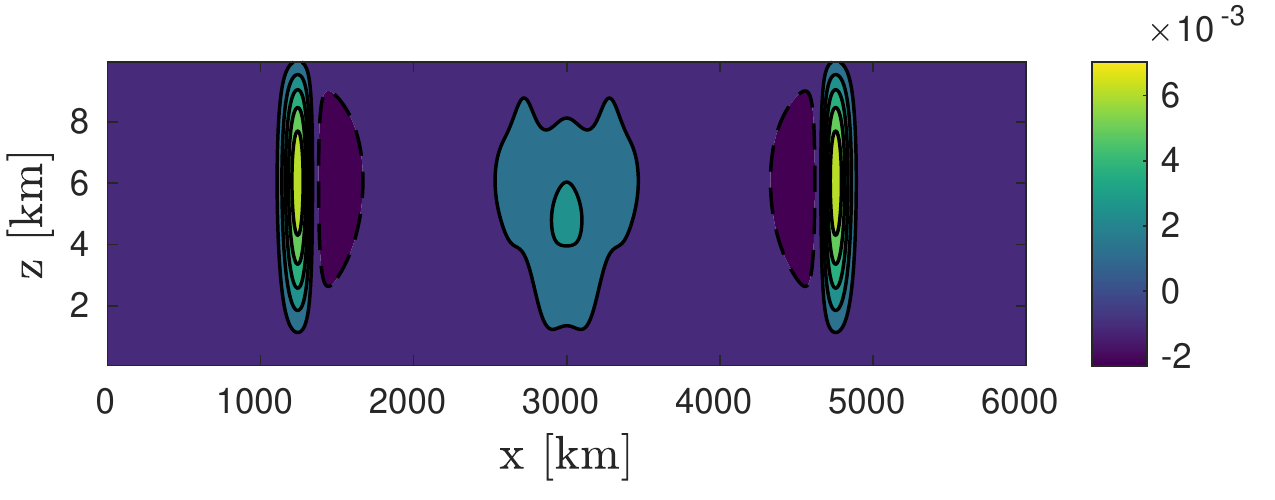}
 \includegraphics[width=.75\columnwidth]{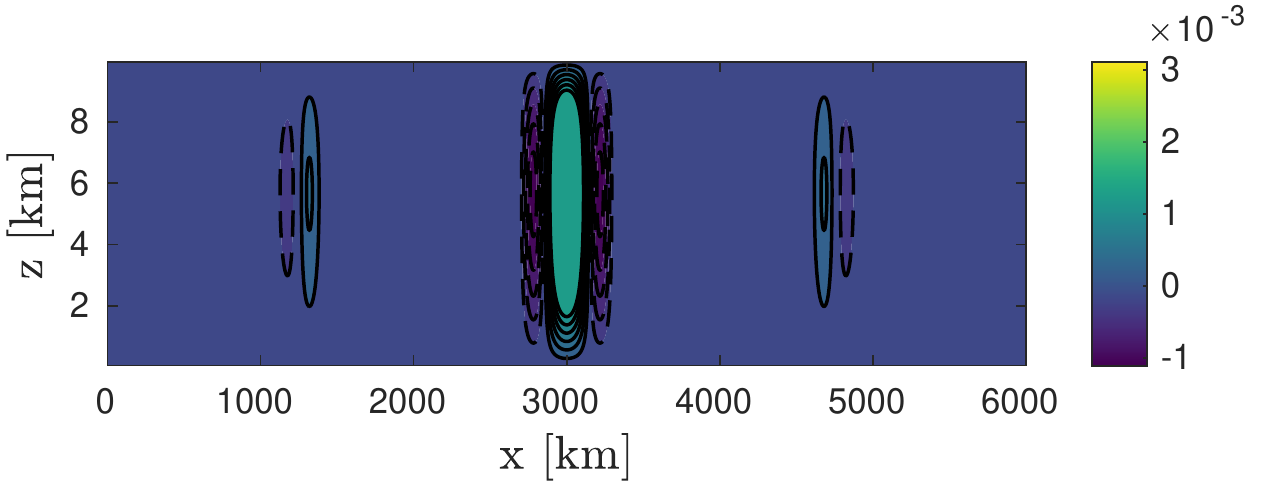}
 \includegraphics[width=.75\columnwidth]{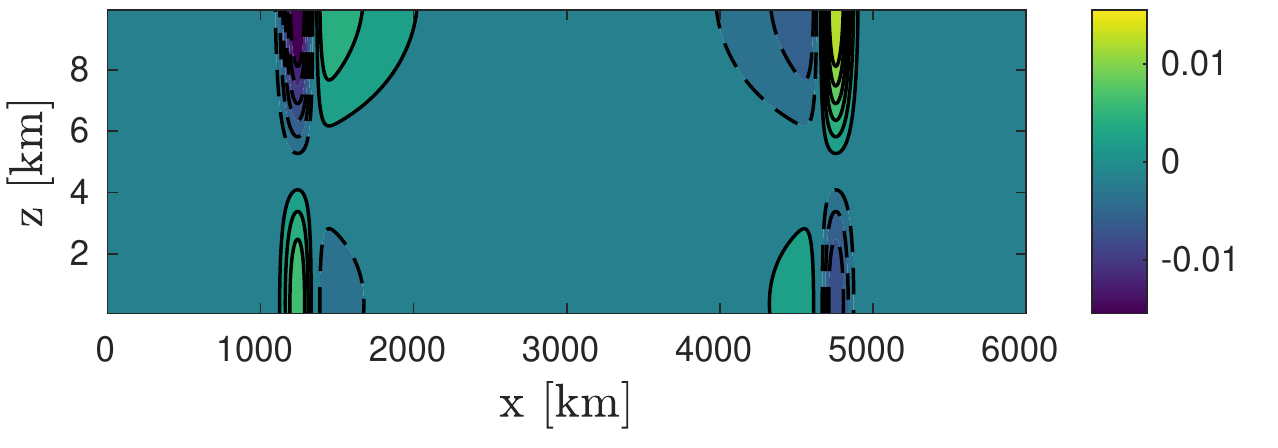}
 \caption{Temperature perturbation (top), vertical velocity (middle) and horizontal velocity (bottom) at final time $T=28800\,\textrm{s}$ for the inertia-gravity wave test with rotation of \cite{BaldaufBrdar2013}, $\Delta x=5\,\textrm{km},\,\Delta z=125\,\textrm{m}$, CFL$_\textrm{adv}=0.9$. Initial perturbation as in Figure \ref{fig:SK94_H},top. Contours in the range $[-6, 6]*10^{-3}\,\textrm{K}$ with a $1.2*10^{-3}\,\textrm{K}$ interval (top), $[-1.2, 1.2]*10^{-3}\,\textrm{ms$^{-1}$}$  with a $2*10^{-4}\,\textrm{ms$^{-1}$}$ interval (middle), $[-0.012, 0.012]\,\textrm{ms$^{-1}$}$ with a $2*10^{-3}\,\textrm{ms$^{-1}$}$ interval (bottom). Negative contours are dashed, zero contours not shown.}
  \label{fig:BB13}
\end{figure}

The initial data trigger a rapidly oscillating vertical acoustic gravity wave pulse that is followed over more than 230 thousand time steps without decay and with small horizontal spread. Superimposed is a longer wavelength internal wave mode that sends two pulses sideways from the center of the initial perturbation, leaving the oscillating acoustic gravity mode behind. Results with the buoyancy-implicit model display good symmetry (Figure \ref{fig:BB13}) and compare well with the reference (Figure 4 in \cite{BaldaufBrdar2013}). The multiscale nature of the case is evident in particular in the plot of the vertical velocity.

\section{Discussion and conclusion}
\label{sec:Conclusions}	

This paper extended a semi-implicit numerical model for the simulation of atmospheric flows to a scheme with time step unconstrained by the internal wave speed and without subtraction of a background state from the primary prognostic variables. The conservative, second-order accurate finite volume discretisation of the rotating compressible equations evolves cell-centered variables through a three-stage procedure, made of an implicit midpoint rule step, an advection step, and an implicit trapezoidal step. By design the model agrees with the pseudo-incompressible system in the small-scale vanishing Mach number limit and with the hydrostatic system at the large-scale limit. Moreover, the discretization is designed so it can be switched straightforwardly to strictly solving either of these two limiting systems. In particular, the hydrostatic mode can be activated simply by flipping a switch in a single equation of the implicit step.

The compressible scheme was applied to a suite of benchmarks of atmospheric dynamics at different scales. Compared with the previous variant of the model in \cite{Benacchio2014, BenacchioEtAl2014}, which used a buoyancy-explicit discretization, the present scheme achieves comparable accuracy, competitive solution quality, and absence of oscillations with much larger time steps for the cases under gravity. New compressible simulations of the hydrostatic-scale inertia-gravity wave tests of \cite{SkamarockKlemp1994} demonstrated the large time step capability of the buoyancy-implicit numerical scheme. An more challenging planetary-scale version of this class of tests was introduced in this paper and revealed the robustness of the discretization for two-hour long time steps. The authors are unaware of published attempts to run the test at this scale. 

An additional test \cite{BaldaufBrdar2013}, geared towards revealing the long-time simulation stability and energy perservation of the scheme, yielded results comparable to those obtained with the reference's higher-order discontinuous Galerkin scheme, albeit with somewhat less of a spreading of the oscillatory mode. The results with the present scheme are superior to those generated by the dynamical core of a weather forecast production code also tested in their paper. 

Furthermore, the hydrostatic- and planetary-scale configurations were run both in pseudo-incompressible mode and in hydrostatic mode, thereby extending the switching capability previously shown in \cite{BenacchioEtAl2014} for the pseudo-incompressible--to--compressible configurations. With increasingly large scales, differences with the compressible runs increased for the pseudo-incompressible runs and decreased for the hydrostatic runs as expected.

The results presented here suggest several avenues of development in a number of areas. First, the scheme serves as the starting point for implementing the multimodel theoretical framework of \cite{KleinBenacchio2016}, which aims to achieve balanced initialization and data assimilation at all scales by smoothly blending between different operation modes. As proposed in \cite{BenacchioEtAl2014}, such a multimodel discretization could be run with reduced soundproof or hydrostatic dynamics during the first time steps after setup or assimilation, then resorting to the fully compressible model for the transient sections. The development in the present work yields hydrostasy at large scale as well as pseudo-incompressibility at small scales as the accessible asymptotic dynamics in the blended scheme. The discretization could then be applied to run tests in spherical geometry, with the ultimate aim of comparing with existing schemes used in numerical weather prediction research and operations.

\section*{Acknowledgments}

T.B. acknowledges funding by the ESCAPE-2 project, European Union's Horizon 2020 research and innovation programme under grant agreement No 800897. R.K.\ acknowledges funding by Deutsche Forschungsgemeinschaft through the Collaborative Research Center CRC~1114 ``Scaling cascades in complex systems'', project~A02, and the support of the European Centre for Medium Range Weather Forecasts under their ECMWF Fellow Program. Extensive discussions with Piotr Smolarkiewicz, Christian K\"uhnlein, and Nils Wedi have been crucial for the present developments. 

\section*{Note}

This work has not yet been peer-reviewed and is provided by the contributing authors as a means to ensure timely dissemination of scholarly and technical work on a noncommercial basis. Copyright and all rights therein are maintained by the authors or by other copyright owners. It is understood that all persons copying this information will adhere to the terms and constraints invoked by each author's copyright. This work may not be reposted without explicit permission of the copyright owner.


\bibliographystyle{plain}
\bibliography{Bibliography}

\end{document}